\documentclass[reqno]{amsart}

\usepackage[normalem]{ulem}

\usepackage{amsmath,amssymb,color}
\usepackage{amsfonts, amscd, epsfig, amsmath, amssymb,enumerate, bbm}
\usepackage{graphicx}
\usepackage{graphics}
\usepackage{color}
\usepackage{mathrsfs}
\usepackage{todonotes}
\usepackage{soul}
\usepackage{faktor}
\numberwithin{equation}{section}

 \usepackage[space]{grffile} 
 \usepackage{etoolbox} 
\usepackage[margin=3cm]{geometry}

\newtheorem{theorem}{Theorem}[section]
\newtheorem{lemma}[theorem]{Lemma}
\newtheorem{proposition}[theorem]{Proposition}

\newtheorem{remark}[theorem]{Remark}

\usepackage{tikz}
\usetikzlibrary{backgrounds}
\usetikzlibrary{patterns,fadings}
\usetikzlibrary{arrows,decorations.pathmorphing}
\usetikzlibrary{decorations}
\usetikzlibrary{calc}
\usetikzlibrary{shapes.misc}

\usepackage{setspace}

\definecolor{light-gray}{gray}{0.95}
\usepackage{float}
\usepackage[colorlinks=true,linkcolor=blue,citecolor=magenta]{hyperref}
\def\centerarc[#1](#2)(#3:#4:#5){\draw[#1] ($(#2)+({#5*cos(#3)},{#5*sin(#3)})$) arc (#3:#4:#5);}

\newcommand{\bb}[1]{{\mathbb #1}}

\newcommand{\<}{\langle}
\renewcommand{\>}{\rangle}

\renewcommand{\epsilon}{\varepsilon}

\newcommand{\R}{\mathbb R}

\newcommand{\Z}{\mathbb Z}

\newcommand{\Prob}{\mathbb P}

\newcommand{\E}{\mathbb E}

\renewcommand{\tilde}{\widetilde}

\newcommand{\normm}[1]{{\left\vert\kern-0.1ex\left\vert\kern-0.1ex\left\vert\; #1 \; \right\vert\kern-0.1ex\right\vert\kern-0.1ex\right\vert}}

\newcommand{\gene}{\mathscr{L}_N}

\renewcommand{\leq}{\leqslant}
\renewcommand{\geq}{\geqslant}

\allowdisplaybreaks 

\title{Moderate deviations for the facilitated exclusion process in  equilibrium}

\author{Linjie Zhao}
\address{School of Mathematics and Statistics \& Hubei Key Laboratory of Engineering Modeling and Scientific Computing, Huazhong University of Science and Technology, Wuhan 430074, China.}
\email{linjie\_zhao@hust.edu.cn}

\thanks{\textbf{Acknowledgments.}  The author would like to thank Oriane Blondel, Cl\'ement Erignoux and Marielle Simon for their helpful discussions.  The author would also like to thank Hugo Da Cunha for sharing the notes on the proof of logarithmic Sobolev inequalities for the FEP. The project is supported by the National Natural Science Foundation of China with grant numbers 12401168 and 12371142, and the Fundamental Research Funds for the Central Universities in China.}


\begin{document}

\maketitle

\begin{abstract}
We derive the moderate deviation principles for the fluctuation fields of the facilitated exclusion process (FEP) in one dimension when the process starts from its stationary measure, both in the symmetric and asymmetric cases. The main step is to prove a super-exponential version of the Boltzmann-Gibbs principle, which relies on the logarithmic Sobolev inequality for the FEP. 
\end{abstract}

\section{Introduction}

The facilitated exclusion process has been widely studied recently after being introduced in \cite{RossiPastorVespignani00}, indicating the existence of a new universality class of nonequilibrium phase transitions in the presence of a conserving quantity.  In the dynamics, there is at most one particle at each site (the exclusion rule), and a particle can jump if and only if the target site is empty and there is a particle nearby (the facilitated rule). Due to the facilitated rule, there exists some critical density $\rho_c$ below which the particles will stay frozen finally. For example, in dimension one, it is easy to see that $\rho_c = 1/2$. However, less is known in higher dimensions, see \cite{erignoux2024hydrodynamic} and references therein.

As far as we know, rigorous results for the FEP have only been proved in one dimension. Stationary states for the FEP were studied in \cite{ayyer2023stationary,Chen2019,Goldstein2020,goldstein2019exact,goldstein2022stationary}. When the FEP starts from the step initial condition, \cite{BaikBarraquandCorwinToufic16} studied the position of the rightmost particle in the totally asymmetric case.   The field of particle
positions around the jump discontinuity was investigated in \cite{barraquand2023weakly} in the weakly asymmetric case. Mixing times for the FEP were considered in \cite{ayre2024mixing,erignoux2024cutoff,massoulie2024cutoff}. Hydrodynamic limits for the FEP were explored in \cite{BESS20,BES21,erignoux2024mapping} with periodic boundary conditions, and in \cite{da2024hydrodynamic,da2025coupling} with open/closed boundary conditions.  Stationary fluctuations were further considered in \cite{erignoux2024stationary}. 

Based on the work of \cite{erignoux2024stationary}, we continue to study moderate deviations for the fluctuation fields of the FEP.  In \cite{erignoux2024stationary}, the authors studied the fluctuation fields of the FEP, 
\[\mathcal{Y}^N_t (du) = \frac{1}{\sqrt{N}} \sum_{x \in \Z} (\eta_x (t) - \rho) \delta_{x/N} (du).\] 
Here, $\eta_x (t)$ is the number of particles at site $x$ at time $t$, $\rho \in (1/2,1)$ is the particle density and $\delta_{x/N}$ is the Kronecker Delta function. The authors proved that when the FEP starts from its stationary measure, the fluctuation fields converge to a generalized Ornstein-Uhlenbeck process in the symmetric case under diffusive scaling. In the asymmetric case under hyperbolic scaling, the limit turns out to be a translation of the initial field.  We rescale the fluctuation fields of the FEP by $a_N$ where $\sqrt{N} \ll a_N \ll N$, which implies that the filed converges to zero as $N \rightarrow \infty$. We are interested in the decay rate and rate function of the above convergence, which is called moderate deviations in the literature. Moderate deviations from hydrodynamic limits have been explored for a few interacting particle systems since the work of \cite{gao2003moderate}, see \cite{wang2006moderate,xue2024nonequilibrium,xue2021moderate,zhao2024moderate,zhao2024wasep}. We prove moderate deviation principles for the rescaled fluctuation fields of the FEP both in symmetric and asymmetric cases. In the symmetric case, the rate function is essentially quadratic; in the asymmetric case, the deviations come only from the initial randomness and the dynamical rate function is zero if it is finite.  See Theorem \ref{thm: mdp}.

The proof follows from the standard approach: for the upper bound, we investigate an exponential martingale related to the process; for the lower bound, we need to study hydrodynamic limits for a perturbed dynamics, see \cite{klscaling}. The main step is to prove a super-exponential version of the Boltzmann-Gibbs principle (see Proposition \ref{prop: boltzmann-gibbs}), which states that one could replace the non-conserved quantities of the dynamics by its density fluctuation field under the correct time scaling \cite[Chapter 11]{klscaling}. The super-exponential Boltzmann-Gibbs principle was previously derived in \cite{zhao2024wasep} for the weakly asymmetric exclusion process and in \cite{zhao2024moderate} for a reaction diffusion model.  The proof in the above two papers relies on the polynomial structure of the exclusion process.  Since the invariant measures of the FEP are not product measures, which is the main difficulty along the proof, we are not aware of how to 
take advantage of the polynomial structure. Instead, our proof is based on the logarithmic Sobolev inequality and equivalence of ensembles for the FEP.  We believe that our approach can be applied to other interacting particle systems, such as the zero range process, as long as the  logarithmic Sobolev inequality holds.

The article is organized as follows. In Section \ref{sec model}, we introduce  the model and its basic properties. The main results are stated in Theorem \ref{thm: mdp}.  Sections \ref{sec: upper bound} and \ref{sec: lower bound} are devoted to the proof of the symmetric case. Since the proof of the asymmetric case is essentially the same as the symmetric case, we only outline it in Section \ref{sec: asymmetric}.

\subsection{Notation} We use $\mathcal{S}$ to denote the space of Schwartz functions on $\R$, and let $\mathcal{S}^\prime$ be its dual space. For a time horizon $T > 0$, let $\mathcal{C}^{1,\infty}_c ([0,T] \times \R)$ be the space of functions on $[0,T] \times \R$ which are continuously differentiable in the first variable and smooth in the second one.  Denote by $\mathcal{D} ([0,T],\mathcal{S}^\prime)$ the space of $\mathcal{S}^\prime$-valued functionals on $[0,T]$ equipped with the Skorokhod topology. For two positive real sequences $\{a_N\}$ and $\{b_N\}$, denote $a_N \ll b_N$ if $\lim_{N \rightarrow \infty} a_N / b_N = 0$. For any measure $\mu$ on some state space, we use $E_{\mu}$ to denote the expectation with respect to $\mu$.

\section{Model and Results}\label{sec model}

In this section, we introduce the model and results rigorously. The model is defined in Subsection \ref{subsec fep}. In Subsection \ref{subsec: inv measure},  we recall basic properties for the invariant measures of the FEP. Main results are stated in Subsection \ref{subsec mdp}.

\subsection{The Facilitated Exclusion Process}\label{subsec fep}

We consider here the one-dimensional  facilitated exclusion process (FEP) \cite{BESS20, BES21}, which is a continuous-time Markov process on the state space $\Omega :=\{0,1\}^{\Z}$. Elements of $\Omega $ are called \emph{configurations}, and will be denoted by $\eta=(\eta_x)_{x\in \Z}$. The FEP is driven by the infinitesimal generator $\mathscr{L}$, acting on local functions $f: \Omega \to \R$ as 
\begin{equation}
\label{eq:LN}
\mathscr{L} f(\eta)= \sum_{x\in \Z}c_{x,x+1}(\eta)\{f(\eta^{x,x+1})-f(\eta)\}.
\end{equation}
Let $N \geq 1$ be the scaling parameter which will go to infinity at last.  We consider two cases in this article:
\begin{itemize}
	\item the symmetric case, where the generator $\gene = N^2 \mathscr{L}$ is accelerated by $N^2$, and the jump rate is given by
	\begin{equation}\label{c01}
		c_{x,x+1}(\eta)=\eta_{x-1}\eta_x(1-\eta_{x+1})+\eta_{x+2}\eta_{x+1}(1-\eta_x);
	\end{equation}
	\item the asymmetric case, where the generator $\gene = N \mathscr{L}$ is accelerated by $N$, and the jump rate is given by
\begin{equation}\label{c02}
	c_{x,x+1}(\eta)=\eta_{x-1}\eta_x(1-\eta_{x+1}).
\end{equation}
\end{itemize}
The configuration $\eta^{x,x+1}$ denotes the configuration where sites $x$ and $x+1$ have been inverted in $\eta$, namely,
\[\eta^{x,x+1}_y:=
\begin{cases}
\eta_y & \mbox{ if } \quad y \neq x,x+1,\\
\eta_x & \mbox{ if } \quad y=x+1,\\
\eta_{x+1} & \mbox{ if } \quad y=x.
\end{cases}\]

\subsection{Invariant measures of the FEP}\label{subsec: inv measure}
Because of the kinetic constraint, the FEP cannot create pairs of neighboring empty sites, because doing so would require a particle located in between to jump out. For this reason, the equilibrium states for the FEP are supported on its \emph{ergodic component} 
\[\mathscr{E} :=\{\eta\in \Omega:  \eta_x+\eta_{x+1}\geq 1, \;\; \forall x\in \Z\},\]
where empty sites are all surrounded by particles. The FEP's grand-canonical states in infinite volume are parametrized by a supercritical particle density $\rho\in (1/2,1]$. Because they only charge the ergodic component $\mathscr{E}$, such states are not product. Although their finite size marginals can be defined by an explicit formula, it will be more convenient in our case to define their finite volume restriction to any finite segment  of the form $\{x,\dots,y\}$ of $\Z$ as the distribution $\Prob_{\{x,\dots,y\}}$ of the following Markov chain. Consider $(\eta_z)_{x\leq z\leq y}$ on the state space $\{0,1\}$,  with initial state $\eta_x \sim {\rm Bernoulli}(\rho)$ and transition probabilities 
\[\Prob(\eta_{z+1}=1\mid \eta_z=1)=1-\Prob(\eta_{z+1}=0\mid \eta_z=1)=d(\rho):=\frac{2\rho-1}{\rho}\]
and 
\[ \Prob(\eta_{z+1}=1\mid \eta_z=0)=1-\Prob(\eta_{z+1}=0\mid \eta_z=0)=1.\]
We can now define the grand-canonical state $\pi_\rho$ as the unique distribution on $\{0,1\}^\Z$ with finite volume marginals given by $\Prob_{\{x,\dots,y\}}$. 

We claim that the construction above is equivalent to the one presented in \cite[Eq. (6.2), p.690]{BESS20}, namely that for an ergodic configuration $\sigma$ on $\Lambda_\ell:=\{1, \dots, \ell\}$ with $p$ particles, 
\begin{equation}
\label{eq:pialpha}
\pi_\rho(\eta_{|\Lambda_\ell}=\sigma)=(1-\rho)d(\rho)^{2p-\ell+1-\sigma(1)-\sigma(\ell)}(1-d(\rho))^{\ell-1-p}.
\end{equation}
Indeed, with the Markovian construction presented above, since $1-d(\rho)=(1-\rho)/\rho$
\begin{align*}
\pi_\rho(\eta_{|\Lambda_\ell}=\sigma)&=\rho^{\sigma_1}(1-\rho)^{1-\sigma_1} d(\rho)^{\sum_{x=2}^\ell \sigma_{x-1}\sigma_x}(1-d(\rho))^{\sum_{x=2}^\ell \sigma_{x-1}(1-\sigma_x)}\\
&=(1-\rho)d(\rho)^{\sum_{x=2}^\ell \sigma_{x-1}\sigma_x}(1-d(\rho))^{-\sigma_1+\sum_{x=2}^\ell \sigma_{x-1}(1-\sigma_x)}.
\end{align*}
Since for an ergodic configuration, all empty sites are preceded by particles, 
\[\sum_{x=2}^\ell \sigma_{x-1}(1-\sigma_x)=\sum_{x=2}^{\ell}(1-\sigma_x)=\ell-1-(p-\sigma_1).\]
Similarly, the number of fully occupied neighboring sites in $\{1,\dots,\ell\}$ can be computed, as the total number of pairs of neighboring sites, minus the pairs containing an empty site
\[\sum_{x=2}^\ell \sigma_{x-1}\sigma_x=\ell-1-\sum_{x=2}^{\ell}\sigma_{x-1}(1-\sigma_x)-\sum_{x=2}^{\ell}\sigma_{x}(1-\sigma_{x-1})=2p-\ell+1-\sigma(1)-\sigma(\ell)\]
as wanted.

\subsection{Moderate deviations for the fluctuation fields}\label{subsec mdp} We fix $\rho \in (1/2,1)$ and start the FEP from the measure $\pi_\rho$. We are interested in the rescaled density fluctuation field $\mu^N_t \in \mathcal{S}^\prime$, which acts on test functions $H \in \mathcal{S}$ as
\[\mu^N_t (H) = \frac{1}{a_N} \sum_{x\in \Z} [\eta_x (t) - \rho ] H(\tfrac{x}{N}),\]  
where $\{a_N\}_{N \geq 1}$ is a sequence of positive reals such that
\[\sqrt{N} \ll a_N \ll N.\]

Fix a time horizon $T > 0$.  Now, we introduce the moderate deviation rate functions. For any local function $f: \Omega \rightarrow \R$, denote
\[\tilde{f} (\rho) = E_{\pi_\rho} [f].\]
For any $\mu \in \mathcal{D} ([0,T],\mathcal{S}^\prime)$ and any $H: [0,T] \times \R \rightarrow \R $ which is continuously differentiable in the time variable and satisfies that $H(t,\cdot), \partial_t H(t,\cdot) \in \mathcal{S}$ for any  $0 \leq t \leq T$, define the linear functional $\ell_t$ as
\[\ell_t (\mu,H) = \mu_t (H_t) - \mu_0 (H_0) - \int_0^t  \mu_s \big((\partial_s + \tilde{h}^\prime (\rho) \Delta) H_s\big) ds,\]
where  $h(\eta) =  \eta_{-1} \eta_0 + \eta_0 \eta_1 - \eta_{-1} \eta_0 \eta_1$.  Note that by definition 
\[ \tilde{h}^\prime (\rho) = \frac{d}{d\rho} E_{\pi_\rho} [h].\]
For the symmetric case, the dynamical moderate deviation rate function is defined as
\[Q_{\rm dyn}^{\rm sym} (\mu):= \sup_{H \in \mathcal{C}^{1,\infty}_c ([0,T] \times \R)} \Big\{   \ell_T (\mu,H) - A(\rho)  \int_0^T \|\nabla H_t\|_{L^2 (\R)}^2 dt \Big\},\]
where \[A(\rho)=\frac{(1-\rho)(2\rho-1)}{\rho}.\] For the asymmetric case, 
 the dynamical moderate deviation rate function is defined as
\[Q_{\rm dyn}^{\rm asym} (\mu):= \sup_{H \in \mathcal{C}^{1,\infty}_c ([0,T] \times \R)} \Big\{   \mu_T (H_T) - \mu_0 (H_0) - \int_0^t  \mu_s \big((\partial_s + A^{\prime} (\rho) \nabla) H_s\big) ds \Big\}.\]
The initial moderate deviation rate function is defined as
\[ Q_{\rm ini} (\mu):= \sup_{\phi \in \mathcal{C}_c^\infty (\R)} \Big\{   \mu_0 (\phi) - \frac{1}{2} B(\rho)  \|\phi\|_{L^2 (\R)}^2 \Big\},\]
where
\[B(\rho) = (2\rho-1)\rho(1-\rho).\]
Finally, the symmetric and asymmetric  rate functions are defined respectively by
\[Q^{\rm sym} (\mu) = Q_{\rm ini} (\mu) + Q^{\rm sym}_{\rm dyn} (\mu), \quad Q^{\rm asym} (\mu) = Q_{\rm ini} (\mu) + Q^{\rm asym}_{\rm dyn} (\mu). \]

Denote by $\bb{P}^N_\rho$ the probability measure on the path space $\mathcal{D} ([0,T],\Omega)$ induced by the FEP with generator $\gene$ with initial measure $\pi_\rho$. Let $\mathbb{E}^N_\rho$ be the corresponding expectation.

Below are the main results of the article.

\begin{theorem}\label{thm: mdp}
The sequence of processes $\{\mu^N_t, 0 \leq t \leq T\}$ satisfies moderate deviation principles with decay rate $a_N^2 / N$ and with rate function $Q^{\rm sym}$ in the symmetric case, and with rate function $Q^{asym}$ in the asymmetric case under the additional assumption $a_N \gg \sqrt{N} (\log N)^2$. 

Precisely speaking, for any closed set $\mathcal{C} \subset \mathcal{D} ([0,T],\mathcal{S}^\prime)$ and any open set $\mathcal{O} \subset \mathcal{D} ([0,T],\mathcal{S}^\prime)$,
\begin{align*}
	\limsup_{N \rightarrow \infty} \frac{N}{a_N^2} \log \bb{P}^N_{\rho} \Big(\{\mu^N_t, 0 \leq t \leq T\} \in \mathcal{C}\Big) \leq - \inf_{\mu \in \mathcal{C}} Q(\mu),\\
		\liminf_{N \rightarrow \infty} \frac{N}{a_N^2} \log \bb{P}^N_{\rho} \Big(\{\mu^N_t, 0 \leq t \leq T\} \in \mathcal{O}\Big) \geq - \inf_{\mu \in \mathcal{O}} Q(\mu),
\end{align*}
where $Q=Q^{\rm sym}$ (respectively $Q=Q^{\rm asym}$) in the symmetric (respectively asymmetric) case.
\end{theorem}

\begin{remark}
	In the asymmetric case, we need the additional assumption on $a_N$ since there is a logarithmic correction to the bound of the equivalence of ensembles for the FEP, see \eqref{equiv ensemble} below.
\end{remark}

\begin{remark}
	In Section \ref{sec: asymmetric}, we shall show that if $Q_{\rm dyn}^{\rm asym} (\mu) < + \infty$, then it must be zero. Thus, in the asymmetric case, the deviations from the dynamics of the FEP is zero under the hyperbolic scaling. To observe deviations from the dynamics, one needs to speed up the process by $N^{3/2}$. We are not aware of how to prove this under this longer time scaling, and leave it as future work.
\end{remark}

\section{The symmetric case: upper bound}\label{sec: upper bound}

In this section, we prove the upper bound in the symmetric case.  We shorten  $Q=Q^{\rm sym}$ and $Q_{\rm dyn}=Q_{\rm dyn}^{\rm sym}$ in this and the next sections since there is no ambiguity. By Feynman-Kac formula (see \cite[page 334]{kipnis1989hydrodynamics} for example), for any test function $H \in  \mathcal{C}^{1,\infty}_c ([0,T] \times \R)$, 
\begin{equation}\label{exp mart}
	\mathscr{M}^N_t (H) := \exp \Big\{  \frac{a_N^2}{N} \Big[\mu^N_t (H_t) - \mu^N_0 (H_0) - \frac{N}{a_N^2} \int_0^t e^{-\tfrac{a_N^2}{N} \mu^N_s (H_s)} (\partial_s + \mathscr{L}_N) e^{\tfrac{a_N^2}{N} \mu^N_s (H_s)} ds \Big] \Big\} 
\end{equation}
is a mean-one exponential martingale. We want to write the martingale as a functional of the rescaled density fluctuation field. 

By direct calculations, 
\begin{multline}\label{cal 1}
 \frac{N}{a_N^2} e^{-\tfrac{a_N^2}{N} \mu^N_s (H_s)} (\partial_s + \mathscr{L}_N) e^{\tfrac{a_N^2}{N} \mu^N_s (H_s)} 
 =   \mu^N_s ( \partial_s H_s) \\
 + \frac{N}{a_N^2} \sum_{x \in \Z} N^2  c_{x,x+1} (\eta(s)) \Big[ \exp \Big\{ \frac{a_N}{N} \big(\eta_{x+1} (s) - \eta_x (s)\big) \big(H(s,\tfrac{x}{N}) - H(s,\tfrac{x+1}{N})\big)\Big\} - 1 \Big].
\end{multline}
Recall that for any local function $g: \Omega \rightarrow \R$, we denote
$\tilde{g} (\rho) = E_{\pi_\rho} [g].$
Note that the FEP is \emph{gradient} in the sense that
\[c_{x,x+1} (\eta) (\eta_x - \eta_{x+1}) = \tau_x h (\eta) - \tau_{x+1} h (\eta),\]
where
$h(\eta) = \eta_{-1} \eta_0 + \eta_0 \eta_1 - \eta_{-1} \eta_0 \eta_1$. Then, by Taylor's expansion and summation by parts formula,
we can write the second term on the righthand side of \eqref{cal 1} as
\begin{multline}\label{cal 2}
	\frac{1}{a_N} \sum_{x \in \Z}\big( \tau_x h(\eta(s)) - \tilde{h} (\rho)\big) \Delta_N H(s,\tfrac{x}{N}) \\
	+ \frac{1}{2N} \sum_{x \in \Z} c_{x,x+1} (\eta(s)) \big(\eta_{x+1} (s) - \eta_x (s)\big)^2 \big(\nabla_N H(s,\tfrac{x}{N})\big)^2  + o_{N,H} (1).
\end{multline}
Here, $o_{N,H} (1)$ denotes some constant depending on $H$ and vanishes as $N \rightarrow \infty$.  Note that we have recentered  the random variable $\tau_x h$ above since $\sum_{x \in \Z} \Delta_N H(s,\tfrac{x}{N}) = 0$.  

In order to further deal with the above two terms, we introduce the following two results. The first one can be viewed as a super-exponential version of the Boltzmann-Gibbs principle (see \cite[page 292]{klscaling}), and the second one concerns the moderate deviations for the invariant measure $\pi_\rho$. Their proofs are presented in Subsections \ref{subsec: superexp boltz} and \ref{subsec: initial mdp} respectively.

\begin{proposition}\label{prop: boltzmann-gibbs}
In the symmetric case,	for any $H \in \mathcal{S}$, for any $\varepsilon > 0$, and for any local function $g: \Omega \rightarrow \R$,
	\begin{equation*}
		\lim_{N \rightarrow \infty} \frac{N}{a_N^2} \log \bb{P}^N_\rho \Big( \sup_{0 \leq t \leq T} \Big| \int_0^t \frac{1}{a_N} \sum_{x\in \Z}  \tau_x V(g,\eta(s)) H(\tfrac{x}{N})  ds\Big| > \varepsilon \Big) = - \infty,
	\end{equation*}
	where 
	\[V(g,\eta) = g(\eta) - \tilde{g} (\rho) - \tilde{g}^\prime (\rho) (\eta_0 - \rho).\]
\end{proposition}

\begin{lemma}\label{lem: mdp pi rho}
	For any $H \in \mathcal{S}$,
	\begin{equation*}
		\lim_{N \rightarrow \infty} \frac{N}{a_N^2} \log E_{\pi_\rho} \Big[ \exp \Big\{  \frac{a_N}{N} \sum_{x \in \Z} (\eta_x - \rho) H(\tfrac{x}{N}) \Big\}\Big] = \frac{1}{2} B(\rho) \|H\|_{L^2 (\R)}^2,
	\end{equation*}
where $B(\rho) = (2\rho-1) \rho (1-\rho) $
\end{lemma}

\begin{remark}
In proposition \ref{prop: boltzmann-gibbs},	we need the supremum over time inside the probability to prove the exponential tightness of the rescaled density fluctuation field.
\end{remark}

\begin{remark}
Using the basic inequality $e^{|x|} \leq e^{x} + e^{-x}$,  Jensen's inequality and Lemma \ref{lem: mdp pi rho}, it is easy to prove that, for any $H \in \mathcal{S}$, for any sequence $r_N$ such that $r_N \gg a_N$, and for any $\varepsilon > 0$,
\[\limsup_{N \rightarrow \infty} \frac{N}{a_N^2} \log \bb{P}^N_\rho \Big( \sup_{0 \leq t \leq T} \Big| \int_0^t \frac{1}{r_N} \sum_{x \in \Z} (\eta_x (s)- \rho) H(\tfrac{x}{N}) ds\Big| > \varepsilon \Big) = - \infty.\]	
\end{remark}

Since 
\[E_{\pi_\rho} [c_{x,x+1} (\eta_x - \eta_{x+1})^2] = \frac{2(1-\rho)(2\rho-1)}{\rho}=: 2 A (\rho),\]
the last two results permit us to rewrite the exponential martingale as
\begin{multline}\label{exp martingale}
	\mathscr{M}^N_t (H) := \exp \Big\{  \frac{a_N^2}{N} \Big[\mu^N_t (H_t) - \mu^N_0 (H_0) - \int_0^t  \mu^N_s \big((\partial_s + \tilde{h}^\prime (\rho) \Delta) H_s\big) ds \\
	- A(\rho)  \int_0^t \|\nabla H_s\|_{L^2 (\R)}^2 ds + \mathscr{R}^N_t (H) \Big] \Big\},
\end{multline}
where $\mathscr{R}^N_t (H)$ is super-exponentially small uniformly in time: for any $\varepsilon > 0$,
\begin{equation}\label{error}
\lim_{N \rightarrow \infty} \frac{N}{a_N^2} \log \bb{P}^N_\rho \big( \sup_{0 \leq t \leq T} \big| \mathscr{R}^N_t (H) \big| > \varepsilon \big) = - \infty.
\end{equation}

In the rest of this section, we prove the upper bound over compact sets in Subsection \ref{subsec: weak upper bound}, and prove the exponential tightness of the sequence $\{\mu^N_t, 0 \leq t \leq T\}$ in Subsection \ref{subsec: exp tightness}.  This concludes the proof of the upper bound.  As stated before, Subsections \ref{subsec: superexp boltz} and \ref{subsec: initial mdp} are devoted to the proof of Proposition \ref{prop: boltzmann-gibbs} and Lemma \ref{lem: mdp pi rho} respectively.

\subsection{Upper bound over compact sets}\label{subsec: weak upper bound} Let $\mathcal{K} \subset \mathcal{D} ([0,T],\mathcal{S}^\prime)$ be compact.  We denote $\mu^N = \{\mu^N_t, 0 \leq t \leq T\}$ for short.  In this subsection, we prove that 
\begin{equation}\label{upper on compact set}
\limsup_{N \rightarrow \infty}  \frac{N}{a_N^2} \log \bb{P}^N_\rho \Big(\mu^N \in \mathcal{K}\Big) \leq - \inf_{\mu \in \mathcal{K}} Q (\mu).
\end{equation}
Since $a_N \gg \sqrt{N}$, for any positive sequences $\{\alpha_N\}, \{\beta_N\}$,
\[\lim_{N \rightarrow \infty} \frac{N}{a_N^2} \log (\alpha_N+ \beta_N) = \max \Big\{ \lim_{N \rightarrow \infty} \frac{N}{a_N^2} \log \alpha_N, \lim_{N \rightarrow \infty} \frac{N}{a_N^2} \log \beta_N \Big\}.\]
Thus, by \eqref{error}, for any $\varepsilon  > 0$,  the term on the leftside of \eqref{upper on compact set}  equals
\begin{align*}
&\limsup_{N \rightarrow \infty}  \frac{N}{a_N^2} \log \bb{P}^N_\rho \Big(\mu^N \in \mathcal{K}, |\mathcal{R}^N_T (H)| \leq \varepsilon \Big)\\
&= \limsup_{N \rightarrow \infty}  \frac{N}{a_N^2} \log \bb{E}^N_\rho \Big[ \mathscr{M}^N_T (H) \mathscr{M}^N_T (H)^{-1} \mathbf{1} \Big\{ \mu^N \in \mathcal{K}, |\mathcal{R}^N_T (H)| \leq \varepsilon \Big\} \Big].
\end{align*}
Note that on the event $\Big\{ \mu^N \in \mathcal{K}, |\mathcal{R}^N_T (H)| \leq \varepsilon \Big\}$, for any $\phi \in \mathcal{C}^\infty_c (\R)$, $	\mathscr{M}^N_T (H)^{-1} $ is bounded from above by
\[\sup_{\mu \in \mathcal{K}} \exp \Big\{ - \frac{a_N^2}{N} \Big[\ell_T (\mu,H)
- A (\rho) \int_0^T \|\nabla H_s\|_{L^2 (\R)}^2 ds -\varepsilon + \mu_0 (\phi) \Big] \Big\} 
 \exp \Big\{  \frac{a_N^2}{N} \mu^N_0 (\phi) \Big\}.\]
Since $\mathscr{M}^N_t (H)$ is a mean-one exponential martingale, for any $\varepsilon > 0$ and for any $H  \in \mathcal{C}^{1,\infty}_c ([0,T] \times \R) ,\phi \in \mathcal{C}^\infty_c (\R)$,  the term on the leftside of \eqref{upper on compact set} is bounded by
\begin{multline*}
	- \inf_{\mu \in \mathcal{K}} \Big\{ \ell_T (\mu,H)
	- A(\rho)  \int_0^T \|\nabla H_s\|_{L^2 (\R)}^2 ds -\varepsilon + \mu_0 (\phi) \Big\} 
	+ \limsup_{N \rightarrow \infty} \frac{N}{a_N^2} \log \bb{E}^N_\rho \Big[ \exp \Big\{  \frac{a_N^2}{N} \mu^N_0 (\phi) \Big\}\Big]\\
	= - \inf_{\mu \in \mathcal{K}} \Big\{\ell_T (\mu,H)
	- A(\rho)  \int_0^T \|\nabla H_s\|_{L^2 (\R)}^2 ds + \mu_0 (\phi) - \frac{1}{2} B(\rho) \|\phi\|_{L^2 (\R)}^2 \Big\} + \varepsilon.
\end{multline*}
In the last identity, we used Lemma \ref{lem: mdp pi rho}.  To conclude the proof, we first let $\varepsilon \rightarrow 0$, then optimize over $H  \in \mathcal{C}^{1,\infty}_c ([0,T] \times \R)$ and $\phi \in \mathcal{C}^\infty_c (\R)$, and finally exchange the order of the supremum and infimum by using the Minimax theorem, see \cite[page 584] {gao2003moderate} for example. 

\subsection{Exponential tightness}\label{subsec: exp tightness}  To prove exponential tightness of the sequence $\{\mu^N\}_{N \geq 1}$, it suffices to check the following two conditions (see \cite[Lemma 3.2]{gao2003moderate} for example): 
\begin{itemize}
	\item[(i)] for any $H \in \mathcal{S}$,
	\begin{equation}\label{exp tight 1}
	\limsup_{M \rightarrow \infty} \limsup_{N \rightarrow \infty} \frac{N}{a_N^2} \log \bb{P}^N_\rho \Big(\sup_{0 \leq t \leq T} |\mu^N_t (H)| > M\Big) = - \infty;
	\end{equation}
	\item[(ii)] for any $H \in \mathcal{S}$ and for any $\varepsilon > 0$,
	\begin{equation}\label{exp tight 2}
	\limsup_{\delta \rightarrow 0} \limsup_{N \rightarrow \infty} \frac{N}{a_N^2} \log \bb{P}^N_\rho \Big(\sup_{0 \leq |t-s| \leq \delta} |\mu^N_t (H) - \mu^N_s (H)| > \varepsilon\Big) = - \infty.
	\end{equation}
\end{itemize}

Recall \eqref{exp martingale}. It is easy to see that the two terms in the second line of \eqref{exp martingale} satisfy the above two conditions. By Lemma \ref{lem: mdp pi rho}, so does the initial density fluctuation field $\mu^N_0 (H)$. Thus, it remains to check the above two conditions for the following martingale term and the time integral term respectively:
\[\frac{N}{a_N^2} \log \mathscr{M}^N_t (H), \quad \int_0^t \mu^N_s (H) ds,\]
where $H \in \mathcal{S}$.

\subsubsection{The martingale term}  For condition (i), by Markov's inequality,
\begin{align*}
	\frac{N}{a_N^2} \log \bb{P}^N_\rho \Big(\sup_{0 \leq t \leq T} \Big|\frac{N}{a_N^2} \log \mathscr{M}^N_t (H) \Big| > M\Big) \leq - M + 	\frac{N}{a_N^2} \log \bb{E}^N_\rho \Big[ \sup_{0 \leq t \leq T} \exp \big\{ \big|\log \mathscr{M}^N_t (H)\big|\} \Big].
\end{align*}
Note that we can first remove the above absolute value on the righthand side by using the basic inequality $e^{|x|} \leq e^x +e^{-x}$. Then, we only need to prove that there exists some constant $C=C(T,H)$ such that 
\[\limsup_{N \rightarrow \infty} \frac{N}{a_N^2} \log \E^N_\rho \Big[ \sup_{0 \leq t \leq T} \mathscr{M}^N_t (H) \Big] \leq C.\]
This follows immediately by using Doob's inequality and the fact that that 
\[\mathscr{M}^N_T (H)^2 \leq e^{C(T,H)a_N^2/N} \mathscr{M}^N_T (2H).\]
Next, we check condition (ii). Since 
\[\sup_{0 \leq t \leq T} \Big|\frac{N}{a_N^2} \log \frac{\mathscr{M}^N_t (H)}{\mathscr{M}^N_{t-} (H)}\Big| \leq \frac{C(H)}{a_N},\]
for $N$ large enough, we have
\[\Big\{ \sup_{0 \leq |t-s| \leq \delta} \Big|\frac{N}{a_N^2} \log \frac{\mathscr{M}^N_t (H)}{\mathscr{M}^N_{s} (H)}\Big| > \varepsilon\Big\} \subset \bigcup_{k=0}^{[T/\delta]} \Big\{ \sup_{k\delta \leq t \leq (k+1)\delta} \Big|\frac{N}{a_N^2} \log \frac{\mathscr{M}^N_t (H)}{\mathscr{M}^N_{k\delta} (H)}\Big| > \varepsilon/4\Big\}. \]
Since the process is stationary in time, 
\begin{multline*}
\limsup_{N \rightarrow \infty} \frac{N}{a_N^2} \log \bb{P}^N_\rho \Big( \sup_{0 \leq |t-s| \leq \delta} \Big|\frac{N}{a_N^2} \log \frac{\mathscr{M}^N_t (H)}{\mathscr{M}^N_{s} (H)}\Big| > \varepsilon \Big) \\
\leq \limsup_{N \rightarrow \infty} \frac{N}{a_N^2} \log \bb{P}^N_\rho \Big( \sup_{0 \leq t \leq \delta} \Big|\frac{N}{a_N^2} \log \mathscr{M}^N_t (H)\Big| > \varepsilon/4 \Big).
\end{multline*}
As before, by Markov's inequality and Doob's inequality, for any $A > 0$, the last limit is bounded by 
\[- \frac{A \varepsilon}{4} +  \limsup_{N \rightarrow \infty} \frac{C(A)N}{a_N^2} \log \bb{E}^N_\rho \Big[\mathscr{M}^N_\delta (H)^A\Big]. \]
Since
\[\mathscr{M}^N_\delta (H)^A \leq \mathscr{M}^N_\delta (AH) \exp \{C(H,A) \delta a_N^2/N\},\]
we conclude the proof by first letting $\delta \rightarrow 0$ and then $A \rightarrow +\infty$.

\subsubsection{The time integral term} We only prove condition (ii) for the integral term, since condition (i) can be proved in the same way. By Garsia-Rodemich-Rumsey inequality (see \cite[page 182]{klscaling} for example), 
\[\sup_{ s \leq t  \leq s+ \delta} \Big| \int_s^t \mu^N_\tau (H) d\tau \Big| \leq C \delta^{\tfrac{1}{3} - \tfrac{1}{6}} B^{\tfrac{1}{12}},\]
where
\[B= \int_0^T dt\,\int_0^T ds\, \frac{1}{|t-s|^{1/3}} \Big|\int_s^t \mu^N_\tau (H) d \tau\Big|.\]
By Markov's inequality, it suffices to show that there exists some constant $C$ such that
\[\limsup_{N \rightarrow \infty} \frac{N}{a_N^2} \log \E^N_\rho \Big[\exp \Big\{  \frac{a_N^2}{N} B  \Big\}\Big] \leq C.\]
As before, we can first remove the absolute value in the expression of $B$. Then, using  Jensen's inequality, the fact that the process is stationary in time and Lemma \ref{lem: mdp pi rho}, the last limit is bounded by
\[\limsup_{N \rightarrow \infty} \frac{1}{T^2}  \int_0^T dt \int_0^T ds  \frac{N}{a_N^2} \log E_{\pi_\rho} \Big[ \exp \Big\{ \frac{a_NT^2 |t-s|^{2/3}}{N}  \sum_{x \in \Z} (\eta_x - \rho) H(\tfrac{x}{N}) \Big\} \Big] \leq C(H,\rho,T), \]
thus concluding the proof.

\subsection{Super-exponential Boltzmann-Gibbs principle}\label{subsec: superexp boltz}  In this subsection, we prove Proposition \ref{prop: boltzmann-gibbs}.  By using Garsia-Rodemich-Rumsey inequality (\cite[page 182]{klscaling}), we can  first remove the supremum over time inside the probability in Proposition \ref{prop: boltzmann-gibbs},  see \cite[Proof of Lemma 3.1]{zhao2024moderate} for details. Thus, we only need to prove that, for any $0 \leq t \leq T$, 
	\begin{equation}\label{eqn boltzman-gibbs}
	\lim_{N \rightarrow \infty} \frac{N}{a_N^2} \log \bb{P}^N_\rho \Big(\Big| \int_0^t \frac{1}{a_N} \sum_{x\in \Z}  \tau_x V(g,\eta(s)) H(\tfrac{x}{N})  ds\Big| > \varepsilon \Big) = - \infty.
\end{equation}

Throughout the proof, we fix a positive integer $\ell = \ell (N)$ such that
	\begin{equation}\label{condition on l}
	\frac{N}{a_N} \ll \frac{\ell}{(\log \ell)^2},  \quad \ell \leq \sqrt{Na_N}.
	\end{equation}  
Such $\ell$ exists by the assumption on $a_N$.  For example, one can take $\ell = N a_N^{-3/4}$. Assume that the value of the local function $g$ depends only on the values of $\{\eta_x, |x| \leq s_g\}$ for some integer $s_g > 0$. Define 
\[V^\ell (g,\eta) = \frac{1}{2\ell^\prime  +1}\sum_{|y| \leq \ell^\prime}\tau_y g(\eta) - \tilde{g} (\rho) - \tilde{g}^\prime (\rho) (\eta_0^\ell - \rho),\]
where $\ell^\prime = \ell -s_g$, and $\eta_0^\ell$ is the average number of particles of the configuration $\eta$ in the segment $\{-\ell,-\ell+1,\ldots,\ell\}$, 
\[\eta_0^\ell = \frac{1}{2\ell+1} \sum_{|y| \leq \ell} \eta_y.\]
We introduce the integer $\ell^\prime$ so that  the value of $V^\ell$ depends only on those of $\{\eta_x, |x| \leq \ell\}$.  Using the summation by parts formula, it is easy to see that 
\[\Big|\frac{1}{a_N} \sum_{x\in \Z}  \tau_x \Big(V(g,\eta) - V^\ell (g,\eta) \Big) H(\tfrac{x}{N})  \Big| \leq \frac{C(H,g) \ell^2}{Na_N}.\]
Since $\ell \ll \sqrt{Na_N}$ by assumption \eqref{condition on l}, we only need to prove that 
\begin{equation*}
	\lim_{N \rightarrow \infty} \frac{N}{a_N^2} \log \bb{P}^N_\rho \Big( \Big| \int_0^t \frac{1}{a_N} \sum_{x\in \Z}  \tau_x V^\ell (g,\eta(s)) H(\tfrac{x}{N})  ds\Big| > \varepsilon \Big) = - \infty.
\end{equation*}
The last identity follows immediately once we can show that, for any $\varepsilon > 0$,
\begin{align}
	\lim_{N \rightarrow \infty} \frac{N}{a_N^2} \log \bb{P}^N_\rho \Big( \Big| \int_0^t \frac{1}{a_N} \sum_{x\in \Z}  \Big(\frac{1}{2 \ell^\prime +1}\sum_{|y-x| \leq \ell^\prime}\tau_y g(\eta(s)) - \tilde{g} (\eta_x^\ell (s))\Big) H(\tfrac{x}{N})  ds\Big| > \varepsilon \Big) 
= - \infty,\label{superexp 1}\\
	\lim_{N \rightarrow \infty} \frac{N}{a_N^2} \log \bb{P}^N_\rho \Big( \Big| \int_0^t \frac{1}{a_N} \sum_{x\in \Z} \Big( \tilde{g} (\eta_x^\ell (s)) - \tilde{g} (\rho) -  \tilde{g}^\prime (\rho) (\eta^\ell_x (s) - \rho) \Big) H(\tfrac{x}{N})  ds\Big| > \varepsilon \Big) 
= - \infty. \label{superexp 2}
\end{align}

In the rest of this subsection, we prove the last two equations respectively. 

\subsubsection{Proof of  \eqref{superexp 1}} We shall use the equivalence of ensembles for the canonical invariant measure of the FEP, which was only proved for densities away from $1/2$ and $1$, see \cite[Proposition 5.6]{erignoux2024stationary}. Thus, the first step is to cut off densities from the boundaries. To this end, we fix some $\delta > 0$ such that $1/2 < \rho - \delta < \rho+ \delta < 1$. We consider the two cases $|\eta_x^\ell (s) - \rho| \geq \delta$  and $|\eta_x^\ell (s) - \rho| < \delta$ respectively. 

For the case $|\eta_x^\ell (s) - \rho| \geq \delta$,  we need to prove that, for any $\varepsilon > 0$,
\begin{multline}\label{superexp 3}
	\lim_{N \rightarrow \infty} \frac{N}{a_N^2} \log \bb{P}^N_\rho \Big( \Big| \int_0^t \frac{1}{a_N} \sum_{x\in \Z}  \Big(\frac{1}{2\ell^\prime  +1}\sum_{|y-x| \leq \ell^\prime}\tau_y g(\eta(s)) - \tilde{g} (\eta_x^\ell (s))\Big)\\
	\times  \mathbf{1} \{|\eta_x^\ell (s) - \rho| \geq \delta\} H(\tfrac{x}{N})  ds\Big| > \varepsilon \Big) 
	= - \infty.
\end{multline}
Without loss of generality, we can assume that $H \geq 0$. Since $g$ is bounded, the absolute value of the term inside the above time integral is bounded by
\[\frac{C(g)}{a_N} \sum_{x\in \Z} \mathbf{1} \{|\eta_x^\ell (s) - \rho| \geq \delta\} H(\tfrac{x}{N}).\]
Then, by Markov's exponential inequality, it suffices to show that, for any $M > 0$,
\begin{equation}\label{markov 1}
\lim_{N \rightarrow \infty} \frac{N}{a_N^2} \log \E^N_\rho \Big[ \exp \Big\{ \int_0^t \frac{C(g)M a_N}{N} \sum_{x\in \Z} \mathbf{1} \{|\eta_x^\ell (s) - \rho| \geq \delta\} H(\tfrac{x}{N}) ds  \Big\}\Big] = 0.
\end{equation} 
Using Jensen's inequality and the fact that $\pi_\rho$ is an invariant measure of the FEP, we only need to prove that 
\begin{equation}\label{boundary density}
\lim_{N \rightarrow \infty} \frac{N}{a_N^2} \log E_{\pi_\rho} \Big[ \exp \Big\{ \frac{C(g)Mta_N}{N} \sum_{x\in \Z} \mathbf{1} \{|\eta_x^\ell  - \rho| \geq \delta\} H(\tfrac{x}{N})   \Big\}\Big] = 0.
\end{equation}
If $\pi_\rho$ was a product measure, the last expression would be true by standard large deviation estimates. In our case, we need more delicate analysis. We first split the sum $\sum_{x \in \Z}$ as $\sum_{|x| > N^2/a_N} + \sum_{|x| \leq N^2/a_N} $.  Since $H \in \mathcal{S}$, there exists some constant $C(H)$ such that
\[\frac{1}{N} \sum_{|x| > N^2/a_N} H(\tfrac{x}{N}) \leq C(H) \frac{a_N^2}{N^2}.\]
Since $a_N \ll N$,
\begin{align*}
&\lim_{N \rightarrow \infty} \frac{N}{a_N^2} \log E_{\pi_\rho} \Big[ \exp \Big\{ \frac{C(g)Mta_N}{N} \sum_{|x| > N^2/a_N} \mathbf{1} \{|\eta_x^\ell  - \rho| \geq \delta\} H(\tfrac{x}{N})   \Big\}\Big] \\
&\leq \lim_{N \rightarrow \infty} C(H,g,M,t,\rho) \frac{a_N}{N} = 0.
\end{align*}
To deal with the sum $\sum_{|x| \leq N^2/a_N}$, we define, for  $0 \leq i \leq 3 \ell$,  \[A_i = \{x \in \Z: x = i + j (3\ell+1)\; \text{ for some } j \in \Z\}. \] 
Then, by H{\"o}lder's inequality, 
\begin{multline*}
\log E_{\pi_\rho} \Big[ \exp \Big\{ \frac{C(g)Mta_N}{N} \sum_{|x| \leq N^2/a_N} \mathbf{1} \{|\eta_x^\ell  - \rho| \geq \delta\} H(\tfrac{x}{N})   \Big\}\Big] \\ \leq \frac{1}{3\ell+1} \sum_{i=0}^{3\ell} \log E_{\pi_\rho} \Big[ \exp \Big\{ \frac{C(g)(3\ell+1)Mta_N}{N} \sum_{x\in A_i, |x| \leq N^2/a_N} \mathbf{1} \{|\eta_x^\ell  - \rho| \geq \delta\} H(\tfrac{x}{N})   \Big\}\Big].
\end{multline*}
 Since $\pi_\rho$ is translation invariant, the last line is bounded by
\[\log E_{\pi_\rho} \Big[ \exp \Big\{ \frac{C \ell Mta_N}{N} \sum_{x \in A_0, |x| \leq N^2/a_N} \mathbf{1} \{|\eta_x^\ell  - \rho| \geq \delta\}   \Big\}\Big]\]
for some constant $C=C(g,H)$. Next, we shall bound the last line by induction and use the Markov property of $\pi_\rho$ introduced in Subsection \ref{subsec: inv measure}. For $y \in \Z$, let $\mathcal{F}_y$ be the $\sigma$-algebra generated by $\{\eta_x, x \leq y\}$.  For convenience, we label the points in $A_0 \cap \{|x| \leq N^2/a_N\}$ from the left to the right as 
\[x_1 < x_2 < \ldots < x_{n},\]
where $n = n(N,\ell) \leq   N^2/(a_N \ell)$. Note that $x_{j+1} = x_j + (3\ell+1)$ for $1 \leq j \leq n-1$.  Define $y_j = x_j + \ell, z_j = x_j - \ell$. By the Markov property of $\pi_\rho$, 
\begin{multline*}
E_{\pi_\rho} \Big[ \exp \Big\{ \frac{C \ell Mta_N}{N} \sum_{j=1}^{n}  \mathbf{1} \{|\eta_{x_j}^\ell  - \rho| \geq \delta\}    \Big\}\Big] \\
= E_{\pi_\rho} \Big[ \exp \Big\{ \frac{C \ell Mta_N}{N} \sum_{j=1}^{n-1}   \mathbf{1} \{|\eta_{x_j}^\ell  - \rho| \geq \delta\}    \Big\} 
 E_{\eta_{y_{n-1}}} \Big( \exp \Big\{ \frac{C \ell Mta_N}{N}   \mathbf{1} \{|\eta_{x_n}^\ell  - \rho| \geq \delta\}    \Big\}\Big)\Big].
\end{multline*}
To bound the above expectation with respect to $E_{\eta_{y_{n-1}}}$, let $\mu_{z_n}$ be the distribution of $\eta_{z_n}$ under the measure $P_{\eta_{y_{n-1}}}$ and let $\nu^1_\rho$ be the Bernoulli measure on $\{0,1\}$ with mean $\rho$.  By exponential decay of finite state Markov chains, there exists some constant $C=C(\rho)$ such that
\[|\frac{\mu_{z_n} (j)}{\nu^1_\rho (j)} - 1| \leq C e^{-C\ell}, \quad j = 0,1.\]
Note also that if $\eta_{z_n} \sim \nu^1_\rho$, then the markov chain $\{\eta_x, x \geq z_n\}$ is stationary. Thus,  
\begin{align*}
	 &E_{\eta_{y_{n-1}}} \Big( \exp \Big\{ \frac{C \ell Mta_N}{N}   \mathbf{1} \{|\eta_{x_n}^\ell  - \rho| \geq \delta\}    \Big\}\Big) \\
	 &= E_{\pi_\rho}  \Big( \frac{d\mu_{z_n}}{d\nu^1_\rho}\exp \Big\{ \frac{C \ell Mta_N}{N}   \mathbf{1} \{|\eta_{x_n}^\ell  - \rho| \geq \delta\}    \Big\}\Big) \\
	 &\leq C \exp \Big\{ \frac{C \ell Mta_N}{N}  - C \ell \Big\}  + E_{\pi_\rho}  \Big(  \exp \Big\{ \frac{C \ell Mta_N}{N}   \mathbf{1} \{|\eta_{x_n}^\ell  - \rho| \geq \delta\}    \Big\}\Big)
\end{align*}
for some constant $C=C(\rho,g,H)$. By \cite{katz1960exponential},  for any $\varepsilon > 0$,
\begin{equation}\label{estimate pi rho}
	P_{\pi_\rho} \Big(|\eta_{x_n}^\ell  - \rho| \geq \varepsilon \Big) \leq C \varepsilon^{-2} e^{-C\ell \varepsilon^2}.
\end{equation}
This allows us to bound
\begin{align*}
	&E_{\pi_\rho}  \Big(  \exp \Big\{ \frac{C \ell Mta_N}{N}   \mathbf{1} \{|\eta_{x_n}^\ell  - \rho| \geq \delta\}    \Big\}\Big) \\
	\leq& 1 +  \exp \Big\{ \frac{C \ell Mta_N}{N}    \Big\} P_{\pi_\rho} \Big(|\eta_{x_n}^\ell  - \rho| \geq \delta \Big) \\
	\leq& 1 + C \delta^{-2}\exp \Big\{ \frac{C \ell Mta_N}{N}  - C \ell \delta^2\Big\} .
\end{align*}
Since $a_N \ll N$, for $N$ large enough, there exists some constant $C=C(g,H,M,t,\delta, \rho)$ such that
\[E_{\eta_{y_{n-1}}} \Big( \exp \Big\{ \frac{C \ell Mta_N}{N}   \mathbf{1} \{|\eta_{x_n}^\ell  - \rho| \geq \delta\}    \Big\}\Big) \leq 1 + C e^{-C\ell}.\]
By induction, 
\[E_{\pi_\rho} \Big[ \exp \Big\{ \frac{C \ell Mta_N}{N} \sum_{j=1}^{n}  \mathbf{1} \{|\eta_{x_j}^\ell  - \rho| \geq \delta\}    \Big\}\Big] \leq  \Big(1 + C e^{-C\ell}\Big)^n \leq \Big(1 + C e^{-C\ell}\Big)^{N^2/(a_N\ell)}.\]
Finally, we bound 
\[\frac{N}{a_N^2}\log E_{\pi_\rho} \Big[ \exp \Big\{ \frac{C \ell Mta_N}{N} \sum_{x \in A_0, |x| \leq N^2/a_N} \mathbf{1} \{|\eta_x^\ell  - \rho| \geq \delta\}   \Big\}\Big] \leq \frac{CN^3}{a_N^3 \ell} e^{- C \ell}\]
for some constant $C=C(g,H,M,t,\delta,\rho)$, which converges to zero since $\ell \gg N /a_N$.  This concludes the proof of  $\eqref{superexp 3}$. 

It remains to prove \eqref{superexp 1} for densities away from $1/2$ and $1$, that is, 
\begin{multline*}
	\lim_{N \rightarrow \infty} \frac{N}{a_N^2} \log \bb{P}^N_\rho \Big( \Big| \int_0^t \frac{1}{a_N} \sum_{x\in \Z}  \Big(\frac{1}{2 \ell^\prime +1}\sum_{|y-x| \leq \ell^\prime}\tau_y g(\eta(s)) - \tilde{g} (\eta_x^\ell (s))\Big)\\
	\times  \mathbf{1} \{|\eta_x^\ell (s) - \rho| < \delta\} H(\tfrac{x}{N})  ds\Big| > \varepsilon \Big) 
	= - \infty.
\end{multline*}
By Markov's inequality, we only need to prove that for any $M > 0$,
\begin{multline}\label{superexp 4}
	\lim_{N \rightarrow \infty} \frac{N}{a_N^2} \log \bb{E}^N_\rho \Big[ \exp \Big\{ \Big| \int_0^t \frac{a_NM}{N} \sum_{x\in \Z}  \Big(\frac{1}{2\ell^\prime +1}\sum_{|y-x| \leq \ell^\prime }\tau_y g(\eta(s)) - \tilde{g} (\eta_x^\ell (s))\Big)\\
\times  \mathbf{1} \{|\eta_x^\ell (s) - \rho| < \delta\} H(\tfrac{x}{N})  ds\Big|  \Big\} \Big]
= 0.
\end{multline}
Since $e^{|x|} \leq e^x + e^{-x}$, $\log (x+y) \leq \log 2 + \max \{\log x, \log y\}$ and $a_N \gg \sqrt{N}$, we can remove the absolute value inside the above exponential. Then, by Feynman-Kac formula \cite[Lemma 7.2, page 336]{klscaling},   the  expression on the left hand side of \eqref{superexp 4} is bounded by
\begin{multline}\label{feynman kac}
	\frac{Nt}{a_N^2} \sup_{f} \Big\{  \int \frac{a_N M}{N} \sum_{x\in \Z}  \Big(\frac{1}{2\ell^\prime+1}\sum_{|y-x| \leq \ell^\prime}\tau_y g(\eta) - \tilde{g} (\eta_x^\ell )\Big)\\
	\times  \mathbf{1} \{|\eta_x^\ell - \rho| < \delta\}  H(\tfrac{x}{N}) f(\eta) \pi_\rho (d \eta) -  N^2 \mathscr{D} (\sqrt{f}; \pi_\rho)\Big\},
\end{multline}
where the supremum is over all local functions $f \geq 0$ on $\Omega$ such that $E_{\pi_\rho} [f] = 1$, and $\mathscr{D} (f; \pi_\rho)$ is the Dirichlet form of $f$ associated with the generator $\mathscr{L}$ with respect to the measure $\pi_\rho$, 
\begin{align*}
	\mathscr{D} (f; \pi_\rho) &= E_{\pi_\rho} [f (- \mathscr{L}) f] 
	= \frac{1}{2} \int_{\Omega} \sum_{x\in \Z} c_{x,x+1} (\eta) \big( f(\eta^{x,x+1}) - f(\eta)\big)^2 \pi_\rho (d \eta).
\end{align*}

The next step is to localize the integral term inside the  supremum in \eqref{feynman kac}.  Define $f_{x,\ell}$ as the  conditional expectation of $f$ conditioned on the $\sigma$-algebra generated by the occupation variables  $\{\eta_{x+i}: -\ell \leq i \leq \ell\}$ with respect to the measure $\pi_\rho$,
\[f_{x,\ell} := E_{\pi_\rho} [f | \eta_{x-\ell}, \eta_{x-\ell +1}, \ldots, \eta_{x+\ell}].\]
Due to the dynamics of the FEP, we also need to consider boundary values of the configuration outside the interval $\Lambda_{x,\ell} := [x-\ell,x+\ell] \cap \Z$. For any $a,b \in \{0,1\}$ and any positive integer $k$, define
\[\pi_{x,\ell,k}^{a,b} (\cdot) = \pi_\rho (\cdot \,| \eta^\ell_x = k/(2\ell+1), \eta_{x-\ell-1}=a, \eta_{x+\ell+1} = b).\]
To make $f_{x,\ell}$ a $\pi_{x,\ell,k}^{a,b}$-density, define
\[f_{x,\ell,k}^{a,b} = \frac{f_{x,\ell}}{E_{\pi_{x,\ell,k}^{a,b}} [f_{x,\ell}]}.\]
With the above notations,
\begin{multline}\label{feynman kac 2}
\int   \Big(\frac{1}{2 \ell^\prime +1}\sum_{|y-x| \leq \ell^\prime}\tau_y g(\eta) - \tilde{g} (\eta_x^\ell )\Big)
\mathbf{1} \{|\eta_0^\ell - \rho| < \delta\} f (\eta) \pi_\rho (d \eta) \\
= \sum_{k: |k/(2\ell+1) - \rho| < \delta, \atop a,b \in \{0,1\}}  m_{x,\ell,k}^{a,b} (f) E_{\pi_{x,\ell,k}^{a,b}} \Big[\Big(\frac{1}{2 \ell^\prime +1}\sum_{|y-x| \leq \ell^\prime} \tau_y g(\eta) - \tilde{g} (\tfrac{k}{2\ell+1} )\Big) f_{x,\ell,k}^{a,b}\Big],
\end{multline}
where the weight function is given by
\[m_{x,\ell,k}^{a,b} (f)= E_{\pi_\rho} \big[f_{x,\ell},\eta^\ell_x = k/(2\ell+1), \eta_{x-\ell-1}=a, \eta_{x+\ell+1} = b \big].\]
Note that
\begin{equation}\label{m sum}
\sum_{k: |k/(2\ell+1) - \rho| < \delta, \atop a,b \in \{0,1\}}  m_{x,\ell,k}^{a,b} (f) \leq 1.
\end{equation}
Then, we can write the integral term inside the  supremum in \eqref{feynman kac} as
\begin{equation}\label{feynman kac integ}
\frac{a_N M}{N} \sum_{x \in \Z} \sum_{k: |k/(2\ell+1) - \rho| < \delta, \atop a,b \in \{0,1\}} H(\tfrac{x}{N}) m_{x,\ell,k}^{a,b} (f) E_{\pi_{x,\ell,k}^{a,b}} \Big[\Big(\frac{1}{2 \ell^\prime +1}\sum_{|y-x| \leq \ell^\prime} \tau_y g(\eta) - \tilde{g} (\tfrac{k}{2\ell+1} )\Big) f_{x,\ell,k}^{a,b}\Big].
\end{equation}

Finally, we shall bound the above term by using the logarithm Sobolev inequality  and equivalence of ensembles for the FEP. By the entropy inequality \cite[page 338]{klscaling}, for any $\gamma > 0$,
\begin{multline*}
	E_{\pi_{x,\ell,k}^{a,b}} \Big[\Big(\frac{1}{2 \ell^\prime +1}\sum_{|y-x| \leq \ell^\prime} \tau_y g(\eta) - \tilde{g} (\tfrac{k}{2\ell+1} )\Big) f_{x,\ell,k}^{a,b} H(\tfrac{x}{N})\Big] \\
	\leq \frac{H(f_{x,\ell,k}^{a,b}| \pi_{x,\ell,k}^{a,b})}{\gamma} + \frac{1}{\gamma} \log E_{\pi_{x,\ell,k}^{a,b}} \Big[ \exp \Big\{ \gamma H(\tfrac{x}{N}) \Big(\frac{1}{2\ell^\prime +1}\sum_{|y-x| \leq \ell^\prime}\tau_y g(\eta) - \tilde{g} (\tfrac{k}{2\ell+1} )\Big) \Big\}\Big],
\end{multline*}
where, for any probability measure $\mu$ and any $\mu$-density $f$, we denote by $H(f|\mu)$ the relative entropy of the density $f$ with respect to $\mu$,
\[H(f|\mu) = E_{\mu} [f \log f].\]
 For any $\pi_{x,\ell,k}^{a,b}$-density $f$, let $\mathscr{D}_{x,\ell,k}^{a,b} (\sqrt{f};\pi_{x,\ell,k}^{a,b})$ be the Dirichlet form of the FEP restricted on the interval $\Lambda_{x,\ell} $ with boundary conditions $\eta_{x-\ell-1} = a, \eta_{x+\ell+1} = b$ with respect to the measure $\pi_{x,\ell,k}^{a,b}$,
\[\mathscr{D}_{x,\ell,k}^{a,b} (\sqrt{f};\pi_{x,\ell,k}^{a,b}) = \frac{1}{2} \int \sum_{y,y+1\in \Lambda_{x,\ell}} c_{y,y+1} (\eta) \{\sqrt{f(\eta^{y,y+1})} - \sqrt{f(\eta)}\}^2 \pi_{x,\ell,k}^{a,b} (d \eta).\]
By convention,
\begin{align*}
	c_{x-\ell,x-\ell+1} (\eta) &= a \eta_{x-\ell} (1-\eta_{x-\ell+1}) +   \eta_{x-\ell+2} \eta_{x-\ell+1} (1-\eta_{x-\ell}),\\
	c_{x+\ell-1,x+\ell} (\eta) &= \eta_{x+\ell-2} \eta_{x+\ell-1} (1-\eta_{x+\ell}) +   b \eta_{x+\ell} (1-\eta_{x+\ell-1}).
\end{align*}
The logarithm Sobolev inequality states that  there exists some constant $C$ independent of $x,\ell, k, a, b$ such that
\[H(f|\pi_{x,\ell,k}^{a,b}) \leq C \ell^2 \mathscr{D}_{x,\ell,k}^{a,b} (\sqrt{f};\pi_{x,\ell,k}^{a,b}).\]
The above inequality can be proved by using the mapping between the ergodic FEP with the symmetric simple exclusion process (SSEP) and the classical logarithm Sobolev inequality for the SSEP, see \cite{Cunha2025}.  Thus, we bound \eqref{feynman kac integ} by
\begin{multline*}
\frac{a_N M}{N} \sum_{x \in \Z} \sum_{k: |k/(2\ell+1) - \rho| < \delta, \atop a,b \in \{0,1\}}  m_{x,\ell,k}^{a,b} (f) \Big\{ \frac{C \ell^2 \mathscr{D}^{a,b}_{x,\ell,k} (\sqrt{f_{x,\ell,k}^{a,b}}; \pi_{x,\ell,k}^{a,b})}{\gamma} \\
+ \frac{1}{\gamma} \log E_{\pi_{x,\ell,k}^{a,b}} \Big[ \exp \Big\{ \gamma H(\tfrac{x}{N}) \Big(\frac{1}{2\ell^\prime +1}\sum_{|y-x| \leq \ell^\prime}\tau_y g(\eta) - \tilde{g} (\tfrac{k}{2\ell+1} )\Big) \Big\}\Big] \Big\}
\end{multline*}
By   convexity of the Dirichlet form, 
\begin{align*}
	&\sum_{x \in \Z} \sum_{k: |k/(2\ell+1) - \rho| < \delta, \atop a,b \in \{0,1\}} m_{x,\ell,k}^{a,b} (f)  \mathscr{D}^{a,b}_{x,\ell,k} (\sqrt{f_{x,\ell,k}^{a,b}} ;\pi_{x,\ell,k}^{a,b}) \\
	\leq& 	\sum_{x \in \Z,\atop a,b \in \{0,1\}}   \mathscr{D}_{x,\ell}^{a,b} (\sqrt{f_{x,\ell}} ;\pi_{\rho}) 
	\leq \sum_{x \in \Z, \atop a,b \in \{0,1\}}   \mathscr{D}_{x,\ell}^{a,b} (f ;\pi_{\rho}) \leq C \ell \mathscr{D}  (f ;\pi_{\rho}).
\end{align*}
Above, for any $\pi_{\rho}$-density $f$, $\mathscr{D}_{x,\ell}^{a,b} (\sqrt{f};\pi_\rho)$ is the Dirichlet form of the FEP restricted on the interval $\Lambda_{x,\ell} $ with boundary conditions $\eta_{x-\ell-1} = a, \eta_{x+\ell+1} = b$ with respect to the measure $\pi_\rho$,
\[\mathscr{D}_{x,\ell}^{a,b} (\sqrt{f};\pi_\rho) = \frac{1}{2} \int \sum_{y,y+1\in \Lambda_{x,\ell}} c_{y,y+1} (\eta) \{\sqrt{f(\eta^{y,y+1})} - \sqrt{f(\eta)}\}^2 \pi_\rho (d \eta).\]
Therefore, \eqref{feynman kac integ} is bounded by
\begin{multline}\label{eqn 2}
\frac{ a_N M}{N \gamma} \Big\{ C \ell^3 \mathscr{D}  (f ;\pi_{\rho}) \\  
+\sum_{x \in \Z} \sup_{k: |k/(2\ell+1) - \rho| < \delta, \atop a,b \in \{0,1\}} \log E_{\pi_{\ell,k}^{a,b}} \Big[ \exp \Big\{ \gamma H(\tfrac{x}{N}) \Big(\frac{1}{2\ell^\prime +1}\sum_{|y| \leq \ell^\prime}\tau_y g(\eta) - \tilde{g} (\tfrac{k}{2\ell+1} )\Big) \Big\}\Big] \Big\},
\end{multline}
where we used \eqref{m sum} and the spatial translation invariance of $\pi_{\rho}$ and we shorten $\pi_{\ell,k}^{a,b} = \pi_{0,\ell,k}^{a,b}$.  To deal with the logarithm in \eqref{eqn 2}, we first recall the equivalence of ensembles \cite[Proposition 5.6]{erignoux2024stationary}, 
\begin{equation}\label{equiv ensemble}
E_{\pi_{\ell,k}^{a,b}} \Big[  \Big| \frac{1}{2\ell^\prime +1}\sum_{|y| \leq \ell^\prime} \tau_y g(\eta) - \tilde{g} (\tfrac{k}{2\ell+1} )\Big| \Big] \leq \frac{C (\log \ell)^2}{\ell}.
\end{equation}
This permits us to bound the logarithm in \eqref{eqn 2} by
\[ \frac{C \gamma (\log \ell)^2}{\ell} H(x/N)
+ \log E_{\pi_{\ell,k}^{a,b}} \Big[ \exp \Big\{ \frac{\gamma H (x/N)}{2\ell^\prime + 1}  \sum_{|y| \leq \ell^\prime} \Big( \tau_y g(\eta) - E_{\pi_{\ell,k}^{a,b}} [\tau_y g] \Big) \Big\}\Big].\]
Intuitively, $(2\ell^\prime+1)^{-1/2} \sum_{|y| \leq \ell^\prime} \Big( \tau_y g(\eta) - E_{\pi_{\ell,k}^{a,b}} [\tau_y g] \Big)$ is approximated by the normal distribution, and thus the second term in the last expression is bounded by 
$C \gamma^2 H(x/N)^2 / \ell$ since
\[E[e^{\alpha X}] \leq C \alpha^2, \quad \forall \alpha \in \R\]
 if $X$ has standard normal distribution, see \cite[Equation (5.21)]{yau1996logarithmic} for its rigorous proof.
Thus, we bound \eqref{feynman kac integ} by
\begin{equation*}
\frac{ C a_N M}{N \gamma} \Big\{  \ell^3 \mathscr{D}  (f ;\pi_{\rho}) 
+ \frac{N \gamma (\log \ell)^2}{\ell} + \frac{N \gamma^2}{\ell}\Big\}.
\end{equation*}
Taking $\gamma = C a_N M \ell^3/N^3$, \eqref{feynman kac} is bounded by
\[C(t,M,H) \Big(\frac{N(\log \ell)^2}{a_N \ell} + \frac{\ell^2}{N^2}\Big).\]
We conclude the proof by \eqref{condition on l}.

\subsubsection{Proof of \eqref{superexp 2}} By Taylor's expansion, 
\[ \Big| \tilde{g} (\eta_x^\ell (s)) - \tilde{g} (\rho) -  \tilde{g}^\prime (\rho) (\eta^\ell_x (s) - \rho) \Big| \leq \frac{1}{2} \sup_{1/2 \leq \alpha \leq 1}|\tilde{g}^{''} (\alpha)| (\eta^\ell_x (s) - \rho)^2.\]
As before, by Markov's inequality, we only need to prove that for any $M > 0$, 
\begin{equation}\label{markov 2}
	\lim_{N \rightarrow \infty} \frac{N}{a_N^2} \log \E^N_{\rho} \Big[\exp \Big\{ \int_0^t \frac{a_N M}{N} \sum_{x \in \Z} (\eta^\ell_x (s) - \rho)^2 H(\tfrac{x}{N}) ds\Big\}\Big] = 0.
\end{equation}
As in the prove of \eqref{markov 1}, we  split the sum $\sum_{x \in \Z}$ into $\sum_{|x| > N^2/a_N} + \sum_{|x| \leq N^2/a_N}$. The first case can be dealt with as before. For the sum $\sum_{|x| \leq N^2/a_N}$, first using Jensen's inequality, then recalling the definitions of $A_i$'s in the proof of \eqref{markov 1}, next  using H{\"o}lder's inequality, and finally using the Markov property of $\pi_\rho$, we bound 
\begin{align*}
	 &\log \E^N_{\rho} \Big[\exp \Big\{ \int_0^t \frac{a_N M}{N} \sum_{|x| \leq N^2/a_N} (\eta^\ell_x (s) - \rho)^2 H(\tfrac{x}{N}) ds\Big\}\Big] \\
&\leq \frac{1}{3\ell+1} \sum_{i=0}^{3\ell}  \log E_{\pi_\rho} \Big[\exp \Big\{ \frac{a_N (3\ell+1) Mt}{N} \sum_{x \in A_i,|x| \leq N^2/a_N} (\eta^\ell_x  - \rho)^2 H(\tfrac{x}{N}) \Big\}\Big]\\
&\leq \frac{C}{\ell}  \sum_{|x| \leq N^2/a_N} \log \Big(\exp \Big\{ \frac{C a_N \ell Mt}{N}   - C \ell\Big\} 
+ E_{\pi_\rho} \Big[\exp \Big\{ \frac{ C a_N \ell Mt}{N}  (\eta^\ell_0  - \rho)^2 H(\tfrac{x}{N})\Big\}\Big] \Big)
\end{align*}
for some constant $C=C(H,\rho)$, see the proof of \eqref{boundary density} for details. 

To bound the expectation inside the above logarithm, we first recall the definition of Sub-Gaussian random variables. We say that a random variable $X$ is  Sub-Gaussian of order $\sigma^2$ if  for any $\theta \in \R$,
\[\log E[e^{\theta X}] \leq \frac{1}{2}\sigma^2 \theta^2.\] 
If $X$ is  Sub-Gaussian of order $\sigma^2$, then for any $0 < \theta < \frac{1}{4\sigma^2}$,
\[ E [e^{\theta X^2}] \leq 1+ 4 \theta \sigma^2,\]
see \cite[Equation (3.14)]{zhao2024wasep}. We claim that $\eta^\ell_0  - \rho$ is Sub-Gaussian of order $C(\rho) \ell^{-1}$. Then, since $a_N \ll N$, for $N$ large enough,
\begin{align*}
	E_{\pi_\rho} \Big[\exp \Big\{ \frac{ C a_N \ell Mt}{N}  (\eta^\ell_0  - \rho)^2  H(\tfrac{x}{N})\Big\}\Big]   \leq 1 + \frac{Ca_N}{N} H(\tfrac{x}{N})
\end{align*}
for some $C = C(H,M,\rho,t)$. Therefore, for $N$ large enough, 
\begin{multline*}
\frac{N}{a_N^2}\log \E^N_{\rho} \Big[\exp \Big\{ \int_0^t \frac{a_N M}{N} \sum_{|x| \leq N^2/a_N} (\eta^\ell_x (s) - \rho)^2 H(\tfrac{x}{N}) ds\Big\}\Big]\\
\leq C  \Big(\frac{N^3}{a_N^3 \ell} e^{-C \ell} + \frac{1}{a_N\ell}\sum_{|x| \leq N^2/a_N} H(\tfrac{x}{N})\Big) \leq C \Big(\frac{N^3}{a_N^3 \ell} e^{-C \ell} + \frac{N}{a_N\ell}\Big),
\end{multline*}
which converges to zero since $\ell \gg N a_N^{-1}$ by \eqref{condition on l}.

It remains to prove the claim. Define 
\[\zeta_x := \rho(\eta_x - \rho) + (1-\rho) (\eta_{x-1} - \rho), \quad x \in \Z.\]
Then, one can easily check that 
\[\eta_x - \rho = \zeta_x + (1-\rho) (\eta_x - \eta_{x-1}), \quad x \in \Z,\]
and that $\{\zeta_x, x \in \Z\}$ is a martingale difference, that is,
\begin{equation*}
	E_{\pi_\rho} [\zeta_x | \mathcal{F}_{x-1}] = 0, \quad \forall x \in \Z,
\end{equation*}
where $\mathcal{F}_x$ is the $\sigma$-algebra generated by $\{\eta_y, y \leq x\}$. Thus, 
\[\eta^\ell_0 - \rho = \zeta^\ell_0  + \frac{1-\rho}{2\ell+1} (\eta_\ell - \eta_{-\ell}).\]
By Hoeffding's lemma \cite[Lemma 2.2, page 27]{boucheron2003concentration}, $\frac{1-\rho}{2\ell+1} (\eta_\ell - \eta_{-\ell})$ is Sub-Gaussian of order $C(\rho) \ell^{-2}$. By using the Markov property of $\pi_{\rho}$ and Hoeffding's lemma again, for any $\theta \in \R$,
\begin{align*}
	\log E_{\pi_{\rho}} [e^{\theta \zeta^\ell_0}] &= \log E_{\pi_{\rho}} \Big[\exp \Big\{ \frac{\theta}{2\ell+1} \sum_{x=-\ell}^{\ell-1} \zeta_x \Big\} E_{\pi_{\rho}} \big[  \exp \big\{ e^{\theta \zeta_\ell / (2\ell+1)}\big\}| \mathcal{F}_{\ell-1}\big]\Big] \\
	&\leq C(\rho) \ell^{-2}  \log E_{\pi_{\rho}} \Big[\exp \Big\{ \frac{\theta}{2\ell+1} \sum_{x=-\ell}^{\ell-1} \zeta_x \Big\} \Big] \leq C(\rho) \ell^{-1},
\end{align*}
where we used the induction argument as before in the last inequality, thus concluding the proof.

\subsection{MDP for the measure $\pi_\rho$}\label{subsec: initial mdp}  In this subsection, we prove Lemma \ref{lem: mdp pi rho}. We follow the ideas in \cite{gao1996moderate}, where MDP for martingales and empirical measures of Markov chains was proved.  Compared with the previous work, first, due to the weight function $H$, the Markov chain here is inhomogeneous; second, since the summation in the lemma is over the integer $\Z$, a cut-off argument is needed.

Recall we write
	\[\eta_x - \rho = \zeta_x + (1-\rho) (\eta_x - \eta_{x-1}),\]
where $\zeta_x := \rho(\eta_x - \rho) + (1-\rho) (\eta_{x-1} - \rho), x \in \Z, $ is  a martingale difference.	Using the summation by parts formula,
	\[\sum_{x \in \Z} (\eta_x - \rho) H(\tfrac{x}{N}) = \sum_{x \in \Z} \zeta_x H(\tfrac{x}{N}) + \frac{1-\rho}{N} \sum_{x \in \Z} \eta_x \nabla_N H (\tfrac{x}{N}),\]
Note that 
\[E_{\pi_\rho} [\zeta_0^2] = (2\rho-1) \rho (1-\rho) =: B(\rho).\]
Using the basic inequality
\[\big|\log E[e^{X+Y}] - \log E[e^X] \big| \leq \|Y\|_{L^\infty (P)},\]
 it suffices to prove that 
\[\lim_{N \rightarrow \infty} \frac{N}{a_N^2} \log E_{\pi_\rho} \Big[ \exp \Big\{ \frac{a_N}{N}\sum_{x \in \Z}\zeta_x H(\tfrac{x}{N}) \Big\}  \Big] = \frac{1}{2} E_{\pi_\rho} [\zeta_0^2] \|H\|_{L^2 (\R)}^2.\]
As before, we first split the sum $\sum_{x \in \Z}$ as $\sum_{|x| > N^2/a_N} + \sum_{|x| \leq N^2/a_N} $. Since 
\[\Big|\frac{1}{N} \sum_{|x| \geq N^2/a_N} \zeta_x H(\tfrac{x}{N}) \Big| \leq \frac{C(H)a_N^2}{N^2},\]
we only need to prove that
\[\lim_{N \rightarrow \infty} \frac{N}{a_N^2} \log E_{\pi_\rho} \Big[ \exp \Big\{ \frac{a_N}{N}\sum_{|x| \leq N^2/a_N} \zeta_x H(\tfrac{x}{N}) \Big\}  \Big] = \frac{1}{2} E_{\pi_\rho} [\zeta_0^2] \|H\|_{L^2 (\R)}^2.\]

\subsubsection{The upper bound} Fix integers $M,L > 0$, which will go to infinity after $N \rightarrow \infty$. Let 
\[K = \Big[ \frac{2N^2a_N^{-1} +1}{L+M}\Big].\]
We divide
\[\Big[- \frac{N^2}{a_N}, \frac{N^2}{a_N}\Big] \bigcap \Z = \bigcup_{j=0}^{K-1} (A_j \cup B_j) \bigcup R,\]
where, for $j = 0, \ldots, K-1$,
\begin{align*}
A_j &= \Big\{ - \Big[ \frac{N^2}{a_N} \Big] + j (M+L) + i: i = 0, \ldots, M-1  \Big\},\\
B_j &= \Big\{ - \Big[ \frac{N^2}{a_N} \Big] + j (M+L) + i: i = M, \ldots, M+L-1  \Big\},
\end{align*}
and $R$ is the set of the remaining points. Note that $|R| \leq L+M$. Fix $p,p^\prime$ such that $\frac{1}{p} + \frac{1}{p^\prime} = 1$. By H{\"o}lder's inequality, for $N$ large enough, there exists some constant $C=C(L+M, \|H\|_\infty)$ such that
\begin{multline*}
E_{\pi_\rho} \Big[ \exp \Big\{ \frac{a_N}{N}\sum_{|x| \leq N^2/a_N} \zeta_x H(\tfrac{x}{N}) \Big\}  \Big]  \leq C \Big(E_{\pi_\rho} \Big[ \exp \Big\{ \frac{a_Np}{N}\sum_{j=0}^{K-1} \sum_{x\in A_j} \zeta_x H(\tfrac{x}{N}) \Big\}  \Big] \Big)^{\frac{1}{p}}\\
\times \Big(E_{\pi_\rho} \Big[ \exp \Big\{ \frac{a_Np^\prime}{N}\sum_{j=0}^{K-1} \sum_{x\in B_j} \zeta_x H(\tfrac{x}{N}) \Big\}  \Big] \Big)^{\frac{1}{p^\prime}}.
\end{multline*}
Let $x_j$ be the rightmost point in $A_j$, 
\[x_j =  - \Big[ \frac{N^2}{a_N} \Big] + j (M+L) + M-1.\]
Using the basic inequality 
\[e^{x} \leq 1 + x + \frac{x^2}{2} + \frac{x^3}{6} e^{|x|},\]
and the fact that $\zeta_x$ is a martingale difference, for $N$ large enough, there exists some constant $C=C(M,p,\|H\|_\infty)$ such that
\begin{align*}
	&E_{\pi_\rho} \Big[ \exp \Big\{ \frac{a_Np}{N} \sum_{x\in A_{j}} \zeta_x H(\tfrac{x}{N}) \Big\}\Big| \mathcal{F}_{x_{j-1}}  \Big]\\
	 \leq& 1 + \frac{a_N^2p^2}{2N^2}  E_{\pi_\rho} \Big[ \Big(\sum_{x\in A_{j}} \zeta_x H(\tfrac{x}{N})\Big)^2 \Big| \mathcal{F}_{x_{j-1}}  \Big]  + \frac{Ca_N^3}{6N^3} \sum_{x \in A_j} |H(\tfrac{x}{N})|^3\\
	 =& 1 + \frac{a_N^2p^2}{2N^2} \sum_{x\in A_{j}} H(\tfrac{x}{N})^2 E_{\pi_\rho} \big[ \zeta_x ^2 \big| \mathcal{F}_{x_{j-1}}  \big]  + \frac{Ca_N^3}{6N^3} \sum_{x \in A_j} |H(\tfrac{x}{N})|^3.
\end{align*}
Since $x - x_{j-1} \geq L$ for $x \in A_j$, by the spatial translation invariance and exponential decay of correlation functions under $\pi_\rho$ \cite[Lemma 6.5]{BESS20}, 
\begin{equation}\label{eqn exp decay}
	\sup_{j} \sup_{x \in A_j} \Big| E_{\pi_\rho} \big[\zeta_x^2 \big| \mathcal{F}_{x_{j-1}}  \big] - E_{\pi_\rho} \big[\zeta_x^2   \big]  \Big| \leq C(\rho) e^{-C(\rho)L}.
\end{equation}
Thus,
\begin{align*}
	&E_{\pi_\rho} \Big[ \exp \Big\{ \frac{a_Np}{N} \sum_{x\in A_{j}} \zeta_x H(\tfrac{x}{N}) \Big\}\Big| \mathcal{F}_{x_{j-1}}  \Big]\\
	\leq& 1 + \frac{a_N^2p^2}{2N^2} \sum_{x\in A_{j}} H(\tfrac{x}{N})^2 \Big(E_{\pi_\rho} \big[ \zeta_x ^2   \big] + C e^{-CL} \Big)+ \frac{Ca_N^3}{6N^3} \sum_{x \in A_j} |H(\tfrac{x}{N})|^3
\end{align*}
Using the induction argument as in the proof of Proposition \ref{prop: boltzmann-gibbs}, 
\begin{align*}
&E_{\pi_\rho} \Big[ \exp \Big\{ \frac{a_Np}{N}\sum_{j=0}^{K-1} \sum_{x\in A_j} \zeta_x H(\tfrac{x}{N}) \Big\}  \Big] \\
\leq& \prod_{j=0}^{K-1} \Big\{ 1 + \frac{a_N^2p^2}{2N^2} \sum_{x\in A_{j}} H(\tfrac{x}{N})^2 \Big(E_{\pi_\rho} \big[ \zeta_x ^2   \big] + C e^{-CL} \Big) + \frac{Ca_N^3}{6N^3} \sum_{x \in A_j} |H(\tfrac{x}{N})|^3\Big\}.
\end{align*}
Similarly, for $N$ large enough, there exists some constant $C=C(p^\prime)$ such that
\begin{equation}\label{eqn 1}
E_{\pi_\rho} \Big[ \exp \Big\{ \frac{a_Np^\prime}{N}\sum_{j=0}^{K-1} \sum_{x\in B_j} \zeta_x H(\tfrac{x}{N}) \Big\}  \Big] \leq \prod_{j=0}^{K-1} \prod_{x \in B_j} \Big(1+ \frac{Ca_N^2}{2N^2} H(\tfrac{x}{N})^2\Big).
\end{equation}
Therefore, 
\begin{align*}
	&\limsup_{N \rightarrow \infty} \frac{N}{a_N^2} \log E_{\pi_\rho} \Big[ \exp \Big\{ \frac{a_N}{N}\sum_{|x| \leq N^2/a_N} \zeta_x H(\tfrac{x}{N}) \Big\}  \Big] \\
	\leq& \limsup_{N \rightarrow \infty} \frac{p}{2N} \sum_{j=0}^{K-1} \sum_{x\in A_{j}} H(\tfrac{x}{N})^2 \Big(E_{\pi_\rho} \big[ \zeta_x ^2   \big] + C e^{-CL} \Big) + \limsup_{N \rightarrow \infty} \frac{C (p^\prime)}{2N} \sum_{j=0}^{K-1}\sum_{x \in B_j} H(\tfrac{x}{N})^2\\
		\leq& \frac{p}{2} E_{\pi_\rho} [\zeta_0^2]  \|H\|_{L^2 (\R)}^2 +  Ce^{-CL} +  \limsup_{N \rightarrow \infty} \frac{C(p^\prime)}{2N} \sum_{j=0}^{K-1}\sum_{x \in B_j} H(\tfrac{x}{N})^2.
\end{align*}
We conclude the proof of the upper bound by first letting $L \rightarrow \infty, L/M \rightarrow 0$, and then letting $p \rightarrow 1$. 

\subsubsection{The lower bound} Fix $p,p^\prime > 0$ such that $\frac{1}{p} + \frac{1}{p^\prime} = 1$. By H{\"o}lder's inequality, for $N$ large enough, there exists some constant $C=C(L+M,\|H\|_\infty)$ such that
\begin{multline*}
	E_{\pi_\rho} \Big[ \exp \Big\{ \frac{a_N}{pN}\sum_{j=0}^{K-1} \sum_{x\in A_j} \zeta_x H(\tfrac{x}{N}) \Big\}  \Big]  \leq C \Big(E_{\pi_\rho} \Big[ \exp \Big\{ \frac{a_N}{N}\sum_{|x| \leq N^2/a_N} \zeta_x H(\tfrac{x}{N}) \Big\}  \Big] \Big)^{\tfrac{1}{p}} \\ \times \Big(E_{\pi_\rho} \Big[ \exp \Big\{ -\frac{a_Np^\prime}{pN}\sum_{j=0}^{K-1} \sum_{x\in B_j} \zeta_x H(\tfrac{x}{N}) \Big\}  \Big] \Big)^{\frac{1}{p^\prime}}.
\end{multline*}
By \eqref{eqn 1}, for $N$ large enough, there exists some constant $C=C(p,p^\prime)$ such that
\[E_{\pi_\rho} \Big[ \exp \Big\{ -\frac{a_Np^\prime}{pN}\sum_{j=0}^{K-1} \sum_{x\in B_j} \zeta_x H(\tfrac{x}{N}) \Big\}  \Big] \leq \prod_{j=0}^{K-1} \prod_{x \in B_j} \Big(1+ \frac{Ca_N^2}{2N^2} H(\tfrac{x}{N})^2\Big).\]
Since $e^x \geq 1+x+\frac{x^2}{2} + \frac{x^3}{6}$ for any $x \in \R$, we have
\begin{align*}
		&E_{\pi_\rho} \Big[ \exp \Big\{ \frac{a_N}{pN} \sum_{x\in A_{j}} \zeta_x H(\tfrac{x}{N}) \Big\}\Big| \mathcal{F}_{x_{j-1}}  \Big]\\
		\geq& 1 + \frac{a_N^2}{2p^2N^2} E_{\pi_\rho} \Big[ \Big(\sum_{x\in A_{j}} \zeta_x H(\tfrac{x}{N})\Big)^2 \Big| \mathcal{F}_{x_{j-1}}  \Big]  + \frac{a_N^3}{6p^3N^3} \Big(\sum_{x\in A_{j}} \zeta_x H(\tfrac{x}{N})\Big)^3\\
		=& 1 + \frac{a_N^2}{2p^2N^2} \sum_{x\in A_{j}} H(\tfrac{x}{N})^2 E_{\pi_\rho} \big[  \zeta_x^2 \big| \mathcal{F}_{x_{j-1}}  \big]  + \frac{a_N^3}{6p^3N^3} \Big(\sum_{x\in A_{j}} \zeta_x H(\tfrac{x}{N})\Big)^3
\end{align*}
By \eqref{eqn exp decay}, the last expression is bounded from below by
\[1 + \frac{a_N^2}{2p^2N^2} E_{\pi_\rho} [\zeta_0^2] \sum_{x\in A_{j}} H(\tfrac{x}{N})^2  - \frac{Ca_N^2}{2p^2N^2} \sum_{x\in A_{j}} H(\tfrac{x}{N})^2 e^{-CL} + \frac{a_N^3}{6p^3N^3} \Big(\sum_{x\in A_{j}} \zeta_x H(\tfrac{x}{N})\Big)^3.\]
Using the induction argument as before, we have
\begin{multline*}
	E_{\pi_\rho} \Big[ \exp \Big\{ \frac{a_N}{pN}\sum_{j=0}^{K-1} \sum_{x\in A_j} \zeta_x H(\tfrac{x}{N}) \Big\}  \Big] \\
	\geq \prod_{j=0}^{K-1} \Big(1 + \frac{a_N^2}{2p^2N^2} E_{\pi_\rho} [\zeta_0^2] \sum_{x\in A_{j}} H(\tfrac{x}{N})^2  - \frac{Ca_N^2}{2p^2N^2} \sum_{x\in A_{j}} H(\tfrac{x}{N})^2 e^{-CL} + \frac{a_N^3}{6p^3N^3} \Big(\sum_{x\in A_{j}} \zeta_x H(\tfrac{x}{N})\Big)^3\Big).
\end{multline*}
Since $\log (1+x) \geq x /(x+1)$ and $a_N \ll N$, 
\begin{align*}
	&\liminf_{N \rightarrow \infty} \frac{N}{a_N^2} \log E_{\pi_\rho} \Big[ \exp \Big\{ \frac{a_N}{N}\sum_{|x| \leq N^2/a_N} \zeta_x H(\tfrac{x}{N}) \Big\}  \Big] \\
	\geq& \frac{1}{2p} \Big(E_{\pi_\rho} [\zeta_0^2]   - C e^{-CL}\Big)	\limsup_{N \rightarrow \infty} \frac{1}{N} \sum_{j=0}^{N-1} \sum_{x \in A_j}  H(\tfrac{x}{N})^2  - \liminf_{N \rightarrow \infty} \frac{C(p,p^\prime)}{N} \sum_{j=0}^{K-1}\sum_{x \in B_j} H(\tfrac{x}{N})^2.
\end{align*}
Finally, we conclude the proof  by first letting $L \rightarrow \infty, L/M \rightarrow 0$, and then letting $p \rightarrow 1$. 

\section{{The symmetric case: lower bound}}\label{sec: lower bound}

For any $\phi \in \mathcal{C}^{\infty}_c (\R)$, define the probability measure $\pi^N_{\rho,\phi}$  on the ergodic component $\mathscr{E}$ as, for any $x < y$,
\begin{multline*}
\log \frac{d \pi^N_{\rho,\phi}}{d \pi_\rho} (\eta_x,\eta_{x+1},\ldots,\eta_y)= \eta_x \log \frac{\rho+\tfrac{a_N}{N} \phi(\tfrac{x}{N})}{\rho} + (1-\eta_x) \log \frac{1-\rho-\tfrac{a_N}{N} \phi(\tfrac{x}{N})}{1-\rho} \\
+ \sum_{x+1\leq z \leq y} \Big\{ \eta_{z-1} \eta_{z} \log \frac{d(\rho+\tfrac{a_N}{N} \phi(\tfrac{z}{N}))}{d(\rho)} + \eta_{z-1} (1-\eta_{z}) \log \frac{\tilde{d}(\rho+\tfrac{a_N}{N} \phi(\tfrac{z}{N}))}{\tilde{d}(\rho)} \Big\},
\end{multline*}
where $d(\rho ) = (2\rho-1) / \rho, \tilde{d} (\rho) = 1 - d(\rho)$.   In other words, for any $x$, under $\pi^N_{\rho,\phi}$,  $\{\eta_y\}_{y \geq x}$ is an inhomogeneous Markov chain on the state space $\{0,1\}$ with initial state $\eta_x \sim {\rm Bernoulli} (\rho+\tfrac{a_N}{N} \phi(\tfrac{x}{N}))$ and transition probabilities
\begin{align*}
\bb{P} (\eta_{z+1} = 1 | \eta_z = 1) = d (\rho+\tfrac{a_N}{N} \phi(\tfrac{z}{N})) = 1 - \bb{P} (\eta_{z+1} = 0 | \eta_z = 1),\\
\bb{P} (\eta_{z+1} = 1 | \eta_z = 0) = 1 = 1 - \bb{P} (\eta_{z+1} = 0 | \eta_z = 0).
\end{align*}
By reflection, the similar result holds true for $\{\eta_y\}_{y \leq x}$.

 Recall that for two probability measures $\mu$ and $\nu$ on $\Omega$, the relative entropy $H(\mu|\nu)$ of $\mu$ with respect to $\nu$ is defined as 
\[H(\mu|\nu) = E_{\mu} [\log (d \mu / d \nu)]\]
if $\mu \ll \nu$, and $H(\mu|\nu) =  + \infty$ otherwise. By direct calculations, we have the following result, whose proof is presented in Appendix \ref{app: pf ini rel entropy}.

\begin{lemma}\label{lem: initial relative entropy}
For any $\phi \in \mathcal{C}^{\infty}_c (\R)$,
\[\lim_{N \rightarrow \infty} \frac{N}{a_N^2} H (\pi^N_{\rho,\phi} |  \pi_\rho) = \frac{1}{2B(\rho)} \|\phi\|_{L^2(\R)}^2.\]
\end{lemma}

For any $H \in \mathcal{C}^{1,\infty}_c ([0,T] \times \R)$ and any $\phi \in \mathcal{C}^{\infty}_c (\R)$, define the measure $\bb{P}^N_{H,\phi}$ on the path space $\mathcal{D} ([0,T], \Omega)$ as 
\[\frac{d \bb{P}^N_{H,\phi}}{d \bb{P}^N_\rho} = \mathscr{M}^N_T (H) \frac{d \pi^N_{\rho,\phi}}{d \pi_\rho}.\]
Under $\bb{P}^N_{H,\phi}$, $(\eta(t))_{t \geq 0}$ is an inhomogeneous Markov process on the state space $\Omega$ with initial distribution $\pi^N_{\rho,\phi}$, whose infinitesimal generator $\mathscr{L}^N_{H,t}$ acts on local functions $f: \Omega \rightarrow \bb{R}$ as
\[\mathscr{L}^N_{H,t}  f(\eta)=N^2\sum_{x\in \Z}c_{x,x+1}(\eta) \exp \big\{ \tfrac{a_N}{N} (\eta_{x+1} - \eta_{x}) (H_t (\tfrac{x}{N} )- H_t (\tfrac{x+1}{N})) \big\}\{f(\eta^{x,x+1})-f(\eta)\},\]
see \cite[Appendix 1.7]{klscaling}.

\subsection{Hydrodynamic limits}\label{subsec: hydrodynamics}  In this subsection, we study the hydrodynamic limits for the density fluctuation field $\mu^N_t$ under the measure $\bb{P}^N_{H,\phi}$. Below is the main result of this subsection.

\begin{proposition}\label{prop: hydrodynaimcs}
The sequence of measures $\{\mu^N_t, 0 \leq t \leq T\}_{N \geq 1}$ converges in $\bb{P}^N_{H,\phi}$-probability, as $N \rightarrow \infty$, to a deterministic measure $\{\alpha (t,u)du, 0 \leq t \leq T\}$, where $\alpha (t,u)$ is the unique weak solution to the following linear PDE:
\begin{equation}\label{hydro eqn}
 \begin{cases}	\partial_t \alpha (t,u) = \tilde{h}^\prime (\rho) \Delta \alpha (t,u) - 2 A(\rho) \Delta H (t,u), &\quad (t,u) \in (0,T] \times \R,\\
\alpha (0,u) = \phi (u), &\quad u \in \R,
\end{cases}
\end{equation}
where $A(\rho) = \frac{(1-\rho)(2\rho-1)}{\rho}$.
\end{proposition}

In the rest of this subsection, we prove the above proposition.  The proof is standard and thus we only sketch it. We start with Dynkin's formula: for any $G \in \mathcal{S} (\R)$, 
\[M^N_t (G) = \mu^N_t (G_t) - \mu^N_0 (G_0) - \int_0^t (\partial_s + \mathscr{L}^N_{H,s}) \mu^N_s (G_s) ds\]
is a martingale with quadratic variation 
\begin{align*}
\<M^N_\cdot (G)\>_t &= \int_0^t \big\{ \mathscr{L}^N_{H,s} \mu^N_s (G_s)^2 - 2 \mu^N_s (G_s) \mathscr{L}^N_{H,s} \mu^N_s (G_s)  \big\}ds\\
&= \int_0^t \Big\{ \frac{N^2}{a_N^2} \sum_{x \in \Z} c_{x,x+1}(\eta(s)) \exp \big\{ \tfrac{a_N}{N} (\eta_{x+1} (s)- \eta_{x} (s)) (H_s (\tfrac{x}{N} )- H_s (\tfrac{x+1}{N})) \big\} \\
&\qquad \times  (G_s (\tfrac{x}{N} )- G_s (\tfrac{x+1}{N}))^2 \Big\}ds \leq C(H,G) t N/a_N^2.
\end{align*}
Since $a_N \gg \sqrt{N}$, by Doob's inequality,
\[\lim_{N \rightarrow \infty} \sup_{0 \leq t \leq T} M^N_t (G) = 0 \quad \text{in $\bb{P}^N_{H,\phi}$-probability.}\]
By direct calculations and Taylor's expansion, the term $(\partial_s + \mathscr{L}^N_{H,s}) \mu^N_s (G_s) $ equals
\begin{align*}
&\mu^N_s (\partial_s G_s) + \frac{N^2}{a_N} \sum_{x \in \Z}  c_{x,x+1}(\eta(s)) \exp \big\{ \tfrac{a_N}{N} 
 (\eta_{x+1} (s)- \eta_{x} (s)) (H_s (\tfrac{x}{N} )- H_s (\tfrac{x+1}{N})) \big\} \\
&\qquad \times (\eta_{x+1} (s)- \eta_{x} (s)) (G_s (\tfrac{x}{N} )- G_s (\tfrac{x+1}{N}))\\
&= \mu^N_s (\partial_s G_s) + \frac{1}{a_N} \sum_{x \in \Z} [\tau_x h (\eta(s)) - \tilde{h} (\rho)] \Delta_N G_s (\tfrac{x}{N}) \\
&\qquad + \frac{1}{N} \sum_{x \in \Z} c_{x,x+1}(\eta(s)) (\eta_{x+1} (s)- \eta_{x} (s))^2 \nabla_N H_s (\tfrac{x}{N}) \nabla_N G_s (\tfrac{x}{N}) + o_N (1),
\end{align*}
where $\lim_{N \rightarrow \infty} o_N (1) = 0$. The aim is to write the above expression as a functional of the density fluctuation field. To this end, we need the following result.

\begin{lemma}\label{lem: 4.3}
For any subset $E \subset \mathcal{D} ([0,T],\Omega)$ such that 
\[\lim_{N \rightarrow \infty} \frac{N}{a_N^2} \log \bb{P}^N_\rho (E) = - \infty,\]
we have
\[\lim_{N \rightarrow \infty} \bb{P}^N_{H,\phi} (E) = 0.\]
\end{lemma}

The proof of the above result is the same as \cite[Lemma 4.3]{zhao2024wasep}, thus is omitted here.  The main ingredients are: 
\begin{itemize}
	\item[(i)]  the entropy inequality;
	\item[(ii)] by \eqref{exp martingale} and Cauchy-Schwarz inequality,
	\[\lim_{N \rightarrow \infty} \frac{N}{a_N^2} \log \bb{P}^N_{H,0} (E) = - \infty;\]
	\item[(iii)]  by Lemma \ref{lem: initial relative entropy}, we have the following relative entropy bound
	\[H (\bb{P}^N_{H,\phi}|\bb{P}^N_{H,0}) = H(\pi^N_{\rho,\phi} | \pi_\rho) \leq C a_N^2/N.\]
\end{itemize}

By Lemma \ref{lem: 4.3}, Proposition \ref{prop: boltzmann-gibbs} and Lemma \ref{lem: mdp pi rho}, we have 
\begin{align*}
M^N_t (G) &= \mu^N_t (G_t) - \mu^N_0 (G_0) - \int_0^t  \mu^N_s ((\partial_s + \tilde{h}^\prime (\rho) \Delta)G_s) ds  \\
&- 2 A(\rho) \int_0^t \int_{\R}  \nabla G_s (u) \nabla H_s (u)\,du\,ds  + \varepsilon^N_t,
\end{align*}
where
\[\lim_{N \rightarrow \infty} \sup_{0 \leq t \leq T} \varepsilon^N_t = 0 \quad \text{in $\bb{P}^N_{H,\phi}$-probability}.\]
By using Lemma \ref{lem: 4.3} and \eqref{exp tight 1}, \eqref{exp tight 2}, it is easy to see that, under $\bb{P}^N_{H,\phi}$, the sequence of measures $\{\mu^N_t, 0 \leq t \leq T\}_{N \geq 1}$ is tight in the space $\mathcal{D} ([0,T], \mathcal{S}^\prime)$.  Moreover, since the process is stationary under $\bb{P}^N_\rho$, following the proof of \cite[Lemma 1.6 and Remark 1.8, Chapter 5]{klscaling} or \cite[Lemma 4.4]{zhao2024wasep}, one can show that any limit point of $\{\mu^N_t, 0 \leq t \leq T\}$ along some subsequence is absolutely continuous with respect to the Lebesgue measure. Denote the density of the limit by $\alpha (t,u)$. Passing to the limit in the above martingale formula, we have
\begin{multline*}
0 = \int_{\R} \alpha (t,u) G (t,u) du - \int_{\R} \alpha (0,u) G (0,u) du - \int_0^t \int_{\R} \alpha (s,u)  (\partial_s + \tilde{h}^\prime (\rho) \Delta)G_s (u)\,du\, ds  \\
	- 2 A(\rho) \int_0^t \int_{\R}  \nabla G_s (u) \nabla H_s (u)\,du\,ds.
\end{multline*}
This implies that $\alpha (t,u)$ is a weak solution to \eqref{hydro eqn}. For the initial profile, it is easy to see that
\[\lim_{N \rightarrow \infty} \frac{1}{a_N} \sum_{x \in \Z} (\eta_x - \rho) G_0 (\tfrac{x}{N}) = \int_{\R} \phi (u) G_0 (u) du\]
in $\pi^N_{\rho,\phi}$-probability.  We conclude the proof by the uniqueness  of the solution to \eqref{hydro eqn}.

\subsection{Proof of the lower bound} For any $G,H \in \mathcal{C}^{1,\infty}_c ([0,T] \times \R)$, we define the scalar product $[\cdot,\cdot]$ as
\[[H,J] = A (\rho) \int_0^T \int_{\R} \nabla H(t,u) \nabla G (t,u)\,du\,dt.\]
Define the equivalence relation $H \sim G$ if $[H,G] =0$. Let $\mathcal{H}$ be the completion of the space $\mathcal{C}^{1,\infty}_c ([0,T] \times \R) / \sim$ with respect to the inner product $[\cdot,\cdot]$. Note that we can rewrite the dynamical rate function as
\[Q_{\rm dyn} (\mu) = \sup_{H \in \mathcal{C}^{1,\infty}_{c} ([0,T] \times \R)}  \Big\{ \ell_T (\mu,H) - [H,H]\Big\}.\]
By Riesz representation theorem, it is easy to  prove the following properties of the rate function, see \cite[Lemma 5.1]{gao2003moderate} for example.

\begin{lemma}\label{lem: rate function representation}
Assume $Q (\mu) < + \infty$. Then, there exist $\phi \in L^2 (\R)$ and $H \in \mathcal{H}$ such that 
\begin{align*}
\mu_0 (\varphi) =  \<\varphi,\phi\>, \; \forall \varphi \in \mathcal{S}; &\quad Q_{\rm ini} (\mu) = \frac{1}{2B(\rho)} \|\phi\|_{L^2 (\R)}^2;\\
\ell_T (\mu,G) = 2[H,G], \; \forall G \in \mathcal{S}; &\quad Q_{\rm dyn} (\mu) = [H,H],
\end{align*}
where $B(\rho) = (2\rho-1)\rho(1-\rho)$.
\end{lemma}

Now, we are ready to conclude the proof of the lower bound. For any open set $\mathcal{O} \subset \mathcal{D} ([0,T],\mathcal{S}^\prime)$, it suffices to show that for any $\mu \in \mathcal{O}$,
\[\liminf_{N \rightarrow \infty} \frac{N}{a_N^2} \log \bb{P}^N_\rho \Big(\{\mu^N_t, 0 \leq t \leq T\} \in \mathcal{O}\Big) \geq - Q (\mu).\]
If $Q(\mu) = + \infty$, then there is nothing to prove. Now, assume $Q(\mu) < + \infty$. Let $H = H(\mu)\in \mathcal{H}$ and $\phi = \phi (\mu) \in L^2 (\R)$ be identified in Lemma \ref{lem: rate function representation}. By using the approximation procedure as in \cite{gao2003moderate}, we can assume that $H \in \mathcal{C}_c^{1,\infty} ([0,T] \times  \R)$, $\phi \in \mathcal{C}^\infty_c (\R)$ and that $\mu(t,u)$ is absolutely continuous with respect to the Lebesgue measure for $0 \leq t \leq T$. Denote the density by $\alpha (t,u)$. Then, by Lemma \ref{lem: rate function representation}, $\alpha$ is the unique solution to the hydrodynamic equation \eqref{hydro eqn}.  

 By Jensen's inequality, 
 \begin{align*}
 	&\liminf_{N \rightarrow \infty} \frac{N}{a_N^2} \log \bb{P}^N_\rho \Big(\{\mu^N_t, 0 \leq t \leq T\} \in \mathcal{O}\Big) \\
 	&=  	\liminf_{N \rightarrow \infty} \frac{N}{a_N^2} \log \bb{E}^N_{H,\phi} \Big[\mathscr{M}^N_T (H)^{-1} \frac{d \pi_\rho}{d \pi^N_{H,\phi}}\mathbf{1}\big( \{\mu^N_t, 0 \leq t \leq T\} \in \mathcal{O} \big) \Big]\\
 	&\geq \limsup_{N \rightarrow \infty} \frac{N}{a_N^2} \log \bb{P}^N_{H,\phi} \big( \{\mu^N_t, 0 \leq t \leq T\} \in \mathcal{O} \big) \\
 	& + \limsup_{N \rightarrow \infty} \frac{N}{a_N^2}  \bb{E}^N_{H,\phi} \big[ \log \mathscr{M}^N_{T} (H)^{-1} \big|  \{\mu^N_t, 0 \leq t \leq T\} \in \mathcal{O} \big] \\
 	& + \limsup_{N \rightarrow \infty} \frac{N}{a_N^2}  \bb{E}^N_{H,\phi} \big[ \log \frac{d \pi_\rho}{d \pi^N_{H,\phi}} \big|  \{\mu^N_t, 0 \leq t \leq T\} \in \mathcal{O} \big].
 \end{align*}
By Proposition \ref{prop: hydrodynaimcs}, the sequence of measures $\{\mu^N_t, 0 \leq t \leq T\}_{N \geq 1}$ converges in $\bb{P}^N_{H,\phi}$-probability, as $N \rightarrow \infty$, to the measure $\{\alpha (t,u)du, 0 \leq t \leq T\}$. Thus, the first term on the right hand side in the last inequality is zero. By Lemma \ref{lem: 4.3}, Proposition \ref{prop: hydrodynaimcs} and Lemma \ref{lem: rate function representation}, the second term equals
\[- (\ell_T (\mu,H) - [H,H]) = - Q_{\rm dyn} (\mu).\] 
By Lemma \ref{lem: initial relative entropy} and Lemma \ref{lem: rate function representation}, the third term equals $- Q_{\rm ini} (\mu)$. This concludes the proof. 

\section{The asymmetric case}\label{sec: asymmetric}

The proof of the asymmetric case is similar to the symmetric one. Thus, we only outline the main differences.

\subsection{The upper bound} We start from the exponential martingale introduced in \eqref{exp mart}. By direct calculations,
\begin{align*}
	&\frac{N}{a_N^2} e^{-\tfrac{a_N^2}{N} \mu^N_s (H_s)}  \mathscr{L}_N e^{\tfrac{a_N^2}{N} \mu^N_s (H_s)} 
	  \\
	=& \frac{N}{a_N^2} \sum_{x \in \Z} N c_{x,x+1} (\eta(s)) \Big[ \exp \Big\{ \frac{a_N}{N}  \big(H_s(\tfrac{x+1}{N}) - H_s(\tfrac{x}{N})\big)\Big\} - 1 \Big]\\
	=& \frac{1}{a_N} \sum_{x \in \Z} [c_{x,x+1} (\eta(s)) - \tilde{c}_{x,x+1} (\rho)] \nabla_N H_s (\tfrac{x}{N}) + o_N (1).
\end{align*}
Our aim is to close the martingale, that is, to express the martingale as a functional of the fluctuation field. To this end, we need the following analogous result of Proposition \ref{prop: boltzmann-gibbs}.

\begin{proposition}
	In the asymmetric case,	 assume further that 
	\[a_N \gg \sqrt{N} (\log N)^2.\]
	Then, for any $H \in \mathcal{S}$, for any $\varepsilon > 0$, and for any local function $g: \Omega \rightarrow \R$,
	\begin{equation*}
		\lim_{N \rightarrow \infty} \frac{N}{a_N^2} \log \bb{P}^N_\rho \Big( \sup_{0 \leq t \leq T} \Big| \int_0^t \frac{1}{a_N} \sum_{x\in \Z}  \tau_x V(g,\eta(s)) H(\tfrac{x}{N})  ds\Big| > \varepsilon \Big) = - \infty,
	\end{equation*}
	where 
	\[V(g,\eta) = g(\eta) - \tilde{g} (\rho) - \tilde{g}^\prime (\rho) (\eta_0 - \rho).\]
\end{proposition}

The proof of the last result is similar to Proposition \ref{prop: boltzmann-gibbs}. The only difference is that when applying the Feynman-Kac formula in \eqref{feynman kac}, we need to minus $N \mathscr{D} (\sqrt{f};\pi_{\rho})$ instead of $N^2 \mathscr{D} (\sqrt{f};\pi_{\rho})$ since the process is speed up by $N$ in the asymmetric case. As a result, by adding up all the estimates in the proof of Proposition \ref{prop: boltzmann-gibbs}, the final bound is given by some constant multiple of
\[\frac{\ell^2}{N a_N} + \frac{N(\log \ell)^2}{a_N \ell} + \frac{\ell^2}{N} + \frac{N^3}{a_N^3 \ell} e^{-C \ell} + \frac{N}{a_N \ell}.\]
We need the last expression go to zero as $N \rightarrow \infty$. One can take $\ell = \varepsilon \sqrt{N}$ for $\varepsilon > 0$, and then let $N \rightarrow \infty, \varepsilon \rightarrow 0$.

Using the last proposition, the exponential martingale introduced in \eqref{exp mart} equals 
\begin{align*}
 	\mathscr{M}^N_t (H) = \exp \Big\{  \frac{a_N^2}{N} \Big[\mu^N_t (H_t) - \mu^N_0 (H_0) - \int_0^t  \mu^N_s \big((\partial_s + A^{\prime} (\rho) \nabla) H_s\big) ds 
+ \mathscr{R}^N_t (H) \Big] \Big\},
\end{align*}
where $A(\rho) = E_{\pi_{\rho}} [\eta_{x-1}\eta_{x} (1-\eta_{x+1})]= (1-\rho)(2\rho-1)/\rho$ and $\mathscr{R}^N_t (H)$ is super-exponentially small uniformly in time in the sense of \eqref{error}. We then conclude the proof by following the proof of the symmetric case line by line.

\subsection{The lower bound} In the asymmetric case, the generator of the perturbed dynamics is given by, for
local functions $f: \Omega \rightarrow \bb{R}$,
\[\mathscr{L}^N_{H,t}  f(\eta)=N \sum_{x\in \Z}c_{x,x+1}(\eta) \exp \big\{ \tfrac{a_N}{N} (H_t (\tfrac{x+1}{N} )- H_t (\tfrac{x}{N})) \big\}\{f(\eta^{x,x+1})-f(\eta)\},\]
where $c_{x,x+1}(\eta) = \eta_{x-1}\eta_{x} (1-\eta_{x+1})$. The following results are respectively parallel to Proposition \ref{prop: hydrodynaimcs} and Lemma \ref{lem: rate function representation}. Since the proof is straightforward, we omit it here. 

\begin{proposition}
	The sequence of measures $\{\mu^N_t, 0 \leq t \leq T\}_{N \geq 1}$ converges in $\bb{P}^N_{H,\phi}$-probability, as $N \rightarrow \infty$, to a deterministic measure $\{\alpha (t,u)du, 0 \leq t \leq T\}$, where $\alpha (t,u)$ is the unique weak solution to the transport equation:
	\begin{equation*}
		\begin{cases}	\partial_t \alpha (t,u) = - A^\prime (\rho) \nabla  \alpha (t,u), &\quad (t,u) \in (0,T] \times \R,\\
			\alpha (0,u) = \phi (u), &\quad u \in \R.
		\end{cases}
	\end{equation*}
\end{proposition}

\begin{lemma}
	If $Q^{\rm asym}_{\rm dyn} (\mu) < + \infty$, then we must have $Q^{\rm asym}_{\rm dyn} (\mu) = 0$.
\end{lemma}

Finally, the lower bound follows from the above two results and the proof of the symmetric case line by line.

\appendix

\section{Proof of Lemma \ref{lem: initial relative entropy}}\label{app: pf ini rel entropy}

In this section, we prove Lemma \ref{lem: initial relative entropy}. Since
\begin{align*}
E_{\pi^N_{\rho,\phi}} [\eta_{x-1} \eta_x] &= \big(\rho+\tfrac{a_N}{N}\phi(\tfrac{x-1}{N})\big) d (\rho+\tfrac{a_N}{N} \phi(\tfrac{x}{N})),  \\
E_{\pi^N_{\rho,\phi}} [\eta_{x-1} (1-\eta_x)] &= \big(\rho+\tfrac{a_N}{N}\phi(\tfrac{x-1}{N})\big) \tilde{d} (\rho+\tfrac{a_N}{N} \phi(\tfrac{x}{N})),
\end{align*}
we have
\begin{multline*}
	H (\pi^N_{\rho,\phi} |  \pi_\rho) = \sum_{x \in \Z} \Big\{ \big(\rho+\tfrac{a_N}{N}\phi(\tfrac{x-1}{N})\big) d (\rho+\tfrac{a_N}{N} \phi(\tfrac{x}{N})) \log \frac{d(\rho+\tfrac{a_N}{N} \phi(\tfrac{x}{N}))}{d(\rho)} \\
	+ \big(\rho+\tfrac{a_N}{N}\phi(\tfrac{x-1}{N})\big) \tilde{d} (\rho+\tfrac{a_N}{N} \phi(\tfrac{x}{N})) \log \frac{\tilde{d}(\rho+\tfrac{a_N}{N} \phi(\tfrac{x}{N}))}{\tilde{d}(\rho)} \Big\}.
\end{multline*}
Note that the last sum is finite since $\phi$ has compact support. Since $\log \frac{d(\rho+\tfrac{a_N}{N} \phi(\tfrac{x}{N}))}{d(\rho)}$ has order $a_N/N$, by replacing $\phi((x-1)/N)$ by $\phi (x/N)$, we have 
\begin{multline*}
	H (\pi^N_{\rho,\phi} |  \pi_\rho) = \sum_{x \in \Z} \Big\{ \big(2\rho-1+\tfrac{2a_N}{N}\phi(\tfrac{x}{N})\big)  \log \frac{d(\rho+\tfrac{a_N}{N} \phi(\tfrac{x}{N}))}{d(\rho)} \\
	+ \big(1-\rho-\tfrac{a_N}{N}\phi(\tfrac{x}{N})\big) \log \frac{\tilde{d}(\rho+\tfrac{a_N}{N} \phi(\tfrac{x}{N}))}{\tilde{d}(\rho)} \Big\}+ \mathcal{O}_\phi (a_N^2/N^2),
\end{multline*}
where $|\mathcal{O}_\phi (a_N^2/N^2)| \leq C(\phi) a_N^2/N^2$.  By Taylor's expansion, 
\begin{align*}
	&\sum_{x \in \Z} \big(2\rho-1+\tfrac{2a_N}{N}\phi(\tfrac{x}{N})\big)  \log \frac{d(\rho+\tfrac{a_N}{N} \phi(\tfrac{x}{N}))}{d(\rho)} \\
	&= \sum_{x \in \Z} \big(2\rho-1+\tfrac{2a_N}{N}\phi(\tfrac{x}{N})\big)   \Big\{ \frac{d(\rho+\tfrac{a_N}{N} \phi(\tfrac{x}{N})) - d(\rho)}{d(\rho)} - \frac{1}{2} \Big(\frac{d(\rho+\tfrac{a_N}{N} \phi(\tfrac{x}{N})) - d(\rho)}{d(\rho)} \Big)^2 \Big\} + \mathcal{O}_\phi (a_N^3/N^2)\\
	&= \sum_{x \in \Z} \Big\{ \frac{a_Nd^\prime(\rho)(2\rho-1)}{Nd(\rho)} \phi (\tfrac{x}{N}) + \frac{a_N^2}{N^2} \phi (\tfrac{x}{N})^2 \Big(\frac{(2\rho-1)d^{\prime \prime} (\rho)}{2d(\rho)} - \frac{2\rho-1}{2} \Big(\frac{d^\prime(\rho)}{d(\rho)}\Big)^2 + \frac{2d^\prime (\rho)}{d (\rho)}\Big) \Big\}\\
	&\qquad + \mathcal{O}_\phi (a_N^3/N^2).
\end{align*}
By direct calculations, 
\[\frac{d^\prime (\rho)}{d (\rho)} = \frac{1}{\rho(2\rho-1)}, \quad \frac{d^{\prime \prime} (\rho)}{d (\rho)} = - \frac{2}{\rho^2 (2\rho-1)}.\]
Then,
\[\frac{(2\rho-1)d^{\prime \prime} (\rho)}{2d(\rho)} - \frac{2\rho-1}{2} \Big(\frac{d^\prime(\rho)}{d(\rho)}\Big)^2 + \frac{2d^\prime (\rho)}{d (\rho)} = \frac{1}{2\rho^2 (2\rho-1)}.\]
Thus,
\begin{align*}
&\sum_{x \in \Z} \big(2\rho-1+\tfrac{2a_N}{N}\phi(\tfrac{x}{N})\big)  \log \frac{d(\rho+\tfrac{a_N}{N} \phi(\tfrac{x}{N}))}{d(\rho)} \\
 =& \sum_{x \in \Z} \Big\{ \frac{a_N}{N \rho} \phi (\tfrac{x}{N})  + \frac{a_N^2}{2N^2 \rho^2 (2\rho-1)} \phi (\tfrac{x}{N})^2 \Big\} +  \mathcal{O}_\phi (a_N^3/N^2).
\end{align*}
Similarly,
\begin{align*}
&\sum_{x \in \Z} \big(1-\rho-\tfrac{a_N}{N}\phi(\tfrac{x}{N})\big) \log \frac{\tilde{d}(\rho+\tfrac{a_N}{N} \phi(\tfrac{x}{N}))}{\tilde{d}(\rho)} \\
=& \sum_{x \in \Z} \Big\{\frac{a_N\tilde{d}^\prime (\rho) (1-\rho)}{N\tilde{d} (\rho)} \phi (\tfrac{x}{N}) 
 + \frac{a_N^2}{N^2} \phi (\tfrac{x}{N})^2 \Big(\frac{(1-\rho)\tilde{d}^{\prime\prime} (\rho)}{2\tilde{d} (\rho)} -\frac{1-\rho}{2} \Big(\frac{\tilde{d}^\prime (\rho)}{\tilde{d} (\rho)}\Big)^2-\frac{\tilde{d}^\prime (\rho)}{\tilde{d} (\rho)}\Big) \Big\} +  \mathcal{O}_\phi (a_N^3/N^2)\\
=&  \sum_{x \in \Z} \Big\{ -\frac{a_N}{N \rho} \phi (\tfrac{x}{N})  + \frac{a_N^2}{2N^2 \rho^2 (1-\rho)} \phi (\tfrac{x}{N})^2 \Big\} +  \mathcal{O}_\phi (a_N^3/N^2).
\end{align*}
where we used the identities
\[\frac{\tilde{d}^\prime (\rho)}{\tilde{d} (\rho)} = - \frac{1}{\rho(1-\rho)}, \quad  \frac{\tilde{d}^{\prime\prime} (\rho)}{\tilde{d} (\rho)} = \frac{2}{\rho^2 (1-\rho)}.\]
Finally, we conclude the proof by noting that
\[\frac{1}{2\rho-1} + \frac{1}{1-\rho} = \frac{\rho}{(2\rho-1) (1-\rho)}.\]

\bibliographystyle{plain}
\bibliography{bibliography.bib}

\end{document}